\documentclass[11pt]{article}
\usepackage{graphicx}
\usepackage{amsmath,amsthm,amssymb}
\usepackage{xypic}
\usepackage{mathabx}
\usepackage{bbm}
\usepackage{color}
\newtheorem{theorem}{Theorem}[section]
\newtheorem{proposition}[theorem]{Proposition}
\newtheorem{definition}[theorem]{Definition}
\newtheorem{example}[theorem]{Example}
\newtheorem{conjecture}[theorem]{Conjecture}
\newtheorem{remark}[theorem]{Remark}
\newtheorem{corollary}[theorem]{Corollary}
\newtheorem{lemma}[theorem]{Lemma}
\newtheorem{warning}[theorem]{Warning}
\newtheorem{principle}[theorem]{Principle}
\newtheorem{question}[theorem]{Question}

\begin{document}

\title{On the Commutative Algebra of Categories}
\author{John D. Berman}
\maketitle

\begin{abstract}
We discuss what it means for a symmetric monoidal category to be a module over a commutative semiring category. Each of the categories of (1) cartesian monoidal categories, (2) semiadditive categories, and (3) connective spectra can be recovered in this way as categories of modules over a commutative semiring category (or $\infty$-category in the last case). This language provides a simultaneous generalization of the formalism of algebraic theories (operads, PROPs, Lawvere theories) and stable homotopy theory, with essentially a variant of algebraic K-theory bridging between the two.
\end{abstract}

\section{Introduction}
Our primary goal in this paper is to develop a categorification of commutative algebra, as a new toolbox for studying symmetric monoidal categories (which we regard in this language as categorified abelian groups). We will first introduce our main results, inspired by familiar facts from ordinary commutative algebra (1.1), before introducing our motivating examples, which come from the study of algebraic theories (1.2). Many of our results can be regarded as theorems in ordinary category theory, and we have sought to use that language as much as possible in the introduction. However, it is unavoidable that we use the language of $\infty$-categories in the body of the paper, for reasons we explain in 1.3.

\subsection{Main results}
Consider the classical situation of the ring $R=\mathbb{Z}[\frac{1}{2}]$. An abelian group $A$ admits the structure of an $R$-module if and only if the `multiplication by 2' homomorphism $2\colon A\rightarrow A$ is invertible. In this case, admitting an $R$-module structure is a \emph{property} of an abelian group, and not extra data. Equivalently, the unique ring homomorphism $\mathbb{Z}\rightarrow R$ induces a fully faithful functor $$\text{CRing}_{R/}\rightarrow\text{CRing}_{\mathbb{Z}/},$$ which is to say that $\mathbb{Z}\rightarrow\mathbb{Z}[\frac{1}{2}]$ is an \emph{epimorphism} of commutative rings.

Here we come upon an odd phenomenon in the category of commutative rings: epimorphisms are \emph{not} the same as surjections. What is true, however, is that epimorphisms of commutative rings coincide with injections of affine schemes (provided the homomorphism is of finite type \cite{EGA} 17.2.6). We review these ideas in Section 2.3. Thus we can identify
\begin{itemize}
\item properties of abelian groups which are classified by actions of commutative rings;
\item affine scheme `injections' into $\text{Spec}(\mathbb{Z})$.
\end{itemize}

Bousfield and Kan call such rings \emph{solid}, and have classified all of them \cite{Core}. As we might expect knowing the geometry of $\text{Spec}(\mathbb{Z})$, they are all built out of quotients and localizations of $\mathbb{Z}$ in a suitable way (quotients or localizations depending on whether the generic point of $\text{Spec}(\mathbb{Z})$ is included in the subset).

In this paper, we prove that many \emph{properties} of symmetric monoidal categories are classified by the actions of solid semiring categories. We should first explain the term \emph{semiring category}.

Most of the categories that arise frequently can be made symmetric monoidal in two different ways. One symmetric monoidal structure (usually the categorical coproduct) behaves additively, the other (usually closed symmetric monoidal) behaves multiplicatively, and the multiplicative structure distributes over the additive structure. For example, Set has a cartesian product which distributes over disjoint union. Ab (abelian groups) has a tensor product which distributes over direct sum. This language is made precise in \cite{GGN}, which we review in Section 2.2.

In particular, each of Set and Ab have the structure of a semiring category (and in a natural way).

\begin{principle}
\label{Princ}
Set and Ab classify the \emph{properties} of symmetric monoidal categories being cocartesian monoidal (respectively additive), just as $\mathbb{Z}[\frac{1}{2}]$ and $\mathbb{Z}/2$ characterize properties of abelian groups.
\end{principle}

What we have said is not literally true. The problem is that Set and Ab are not finitely generated as semiring categories, and will therefore have poor algebraic properties. However, there are two ways to resolve this problem, and either way the principle becomes true.

The first option is to note that Set and Ab are finitely generated under colimits, which is essentially to say they are \emph{presentable} (or \emph{locally presentable} for some authors). Lurie has developed a commutative algebra of presentable $\infty$-categories (\cite{HA} 4.8) and has proven $\infty$-categorical analogues of the statements of Principle \ref{Princ}. (The 1-categorical statements are also true, and can be proven using the same techniques applied to the $\infty$-category of 1-categories.)

There are many benefits to this first approach. Presentable categories (and $\infty$-categories) have excellent formal properties, guaranteed by results like the Yoneda lemma and the adjoint functor theorem, which make them an ideal setting for proving universal properties. For example, Lurie uses his commutative algebra of presentable $\infty$-categories to define a well-behaved symmetric monoidal smash product of spectra (\cite{HA} 4.8.2). This is in answer to a long-standing open problem of late twentieth century homotopy theory; the first solution \cite{EKMM} predates Lurie's solution by a decade, but his is the first from an $\infty$-categorical perspective and is surprisingly slick compared with its predecessors.

The second approach to Principle \ref{Princ} is to insist on working just with finitely generated semiring categories, and take the smallest subsemirings of Set and Ab; that is, the subcategories generated by sums and products of the additive and multiplicative units. In the case of Set, we recover in this way the semiring category Fin of finite sets. 

In the case of Ab, we recover the subcategory of finitely generated free abelian groups. This is equivalent to the Burnside category Burn, which can also be described via a group-completed span construction applied to Fin (see Example \ref{BurnSub}). We can also apply a span construction without group-completing; this is the effective Burnside category $\text{Burn}^\text{eff}$, and is equivalent to the category of finitely generated free commutative monoids. Our main result is as follows:

\begin{theorem}
\label{IntroThm1}
Each of the following is a semiring category: Fin, $\text{Fin}^\text{op}$, $\text{Fin}^\text{iso}$, $\text{Fin}^\text{inj}$, $\text{Fin}^{\text{inj},\text{op}}$, $\text{Fin}_\ast$, $\text{Fin}_\ast^\text{op}$, $\text{Burn}^\text{eff}$, and Burn. Superscripts denote that we are only allowing injections (inj) or bijections (iso). $\text{Fin}_\ast$ denotes pointed finite sets.

Moreover, for each semiring category $\mathcal{R}$ on this list, the forgetful functor $\text{Mod}_R\rightarrow\text{SymMon}$ is fully faithful; that is, being an $\mathcal{R}$-module is a property of a symmetric monoidal category, rather than containing extra structure. These properties are as follows: \\
\begin{tabular}{l r}
$\mathcal{R}$ &$\text{Mod}_\mathcal{R}$ \\ \hline
$\text{Fin}^\text{iso}$ &arbitrary \\
$\text{Fin}$ &cocartesian monoidal \\
$\text{Fin}^\text{op}$ &cartesian monoidal \\
$\text{Fin}^\text{inj}$ &semi-cocartesian monoidal \\
$\text{Fin}^{\text{inj},\text{op}}$ &semi-cartesian monoidal \\
$\text{Fin}_\ast$ &cocartesian monoidal with zero object \\
$\text{Fin}_\ast^\text{op}$ &cartesian monoidal with zero object \\
$\text{Burn}^\text{eff}$ &semiadditive (both cocartesian and cartesian monoidal) \\
$\text{Burn}$ &additive \\
\end{tabular}
\end{theorem}

However, we will delay the last result (that Burn-modules are additive categories) until a sequel \cite{Berman2}. Although it is possible to prove using similar techniques to those in this paper (and we invite the reader to do so), the proof is more ad hoc. In the sequel, we will develop tools from which this result naturally follows.

We also prove all these results for $\infty$-categories. They remain true verbatim, and there is an additional result as well (which does not have a natural 1-categorical analogue!):

\begin{theorem}
\label{IntroThm2}
There is a (solid) semiring $\infty$-category $\vec{\mathbb{S}}$ for which $\text{Mod}_{\vec{\mathbb{S}}}$ is equivalent to the $\infty$-category of connective spectra.
\end{theorem}

The idea of identifying spectra with $\mathbb{S}$-modules in a larger category should be familiar to stable homotopy theorists; in fact, standard constructions of spectra (like Elmendorf-Kriz-Mandell-May \cite{EKMM}) take this approach.\footnote{This theorem should not itself be taken as a construction or definition of spectra. We cannot work seriously with symmetric monoidal $\infty$-categories unless we already know something about connective spectra, so such a definition would be circular.}

A benefit of such an approach is that it provides for us a way to compare symmetric monoidal $\infty$-categories to spectra, internal to our categorified commutative algebra. For example, we have a free functor from symmetric monoidal $\infty$-categories to spectra, given by tensoring with $\vec{\mathbb{S}}$. This operation can be regarded as a relative of algebraic K-theory, although not quite the same.

Roughly, if $\mathcal{C}$ is a symmetric monoidal $\infty$-category, then $\mathcal{C}\otimes\vec{\mathbb{S}}$ is obtained from $\mathcal{C}$ first by formally inverting all morphisms (taking the classifying space) and then group-completing. In contrast, most constructions of algebraic K-theory (such as \cite{Mandell}) operate by throwing out all non-invertible morphisms and then group-completing.

Various notions of categorified rings have appeared before, including to study iterated K-theory \cite{BDRR}, Tannaka duality \cite{FunProGpoid}, and smash products of spectra (\cite{HA} 4.8). The framework in this paper is general enough to permit comparisons to many of the other notions of categorified rings; for the most part, we will not discuss this, but the comparison to what may be called \emph{presentable} categorified rings (as in \cite{HA} 4.8 and \cite{FunProGpoid}) is touched on in Remark \ref{PK} and Lemma \ref{SMYoneda}, and the author hopes to elaborate on this relationship in a sequel \cite{Berman2}.

\subsection{Applications}
When we wish to work with a particular type of algebraic structure (such as a group, abelian group, ring, Lie algebra, etc.), it is often the case that the axioms of that structure can be encoded most effectively in the data of a symmetric monoidal category.

For example, if we would like to define a commutative monoid object internal to a symmetric monoidal category $\mathcal{C}^\odot$, we may demand that we have an object $X$, a binary operation $X\odot X\rightarrow X$, a unit $1\rightarrow X$, and various axioms. Or we may simply ask for a symmetric monoidal functor $F\colon\text{Fin}^\amalg\rightarrow\mathcal{C}^\odot$. The two are equivalent; in particular, $F$ will send $n$ to $X^{\odot n}$, the map $2\rightarrow 1$ to the binary operation, and $0\rightarrow 1$ to the unit. All other morphisms in Fin are generated by these, and associativity, commutative, and unitality are encoded in Fin itself.

More generally, we can imagine other symmetric monoidal categories $\mathcal{T}$, whose objects are still labeled by finite sets (and the symmetric monoidal operation is still given by disjoint union of finite sets), but that may have other types of morphisms. All of the semiring categories of Theorem \ref{IntroThm1} are of this form. A symmetric monoidal functor $F\colon\mathcal{T}\rightarrow\mathcal{C}^\odot$ will always consist of some object $F(1)\in\mathcal{C}$ with various kinds of units, binary operations, ternary operations, etc. Many kinds of algebraic structure (including the structure of an algebra over any operad -- see Remark \ref{OperadExample}) can be encoded in this way.

The only constraint we have put on $\mathcal{T}$ is that it be \emph{cyclic} -- that it be generated by a single object under its symmetric monoidal operation. Traditionally such a $\mathcal{T}$ is called a PROP -- for \emph{product and permutation category} \cite{PROP}. If $\mathcal{T}$ is also cartesian monoidal, it is called a Lawvere theory \cite{Lawvere}. This formalism of algebraic theories has been used extensively, particularly in homotopy theory. For early references, see \cite{BV} for infinite loop spaces and \cite{SegalPROP} for the Segal PROP.

In light of Theorem \ref{IntroThm1}:

\begin{corollary}
\label{CorIntro}
PROPs are cyclic $\text{Fin}^\text{iso}$-modules, and Lawvere theories are cyclic $\text{Fin}^\text{op}$-modules.
\end{corollary}

For example, the Lawvere theory modeling associative algebras is $\text{Fin}_\ast^\text{op}$, while the Lawvere theory modeling commutative algebras is $\text{Burn}^\text{eff}$.

In general, we identify cyclic modules over various semiring categories with different flavors of algebraic theories, all of which can be folded into the study of categorified commutative algebra.

\begin{example}
\label{IntroEx}
Another example comes from equivariant homotopy theory. Given a space with the action of a finite group $G$, its homotopy groups will naturally inherit the structure of product-preserving functors $$\text{Fin}_G^\text{op}\rightarrow\text{Set}.$$ Given a generalized homology theory (or spectrum) with a $G$-action, its values at a $G$-space will inherit even more structure, namely that of a product-preserving functor $$\text{Burn}_G\rightarrow\text{Ab}.$$ $\text{Burn}_G$ is the category of spans of finite $G$-sets.

And if the generalized homology theory has well-behaved cup products (is a commutative ring spectrum), we will have even more structure than this: that of a symmetric monoidal functor from a bispan category $\text{Poly}_G$.

We can regard $\text{Fin}_G^\text{op}$, $\text{Burn}_G$, and $\text{Poly}_G$ as all being cyclic modules over $\text{Fin}_G^\text{op}$ -- that is, equivariant Lawvere theories. They play an even more central role in equivariant homotopy theory than we have indicated: in fact, by a result of Elmendorf (\cite{Elmendorf} Theorem 1), the $\infty$-category of $G$-equivariant spaces is equivalent to space-valued models of the $G$-Lawvere theory $\text{Fin}_G^\text{op}$. And by a result of Guillou-May (\cite{GMay2} Theorem 0.1, but see also \cite{BarMack1} Example B.6 for this language) the $\infty$-category of $G$-equivariant spectra is equivalent to spectra-valued models of the $G$-Lawvere theory $\text{Burn}_G$:\footnote{These results of Elmendorf and Guillou-May are not only theoretical. Most examples of equivariant spaces and spectra which are algebraic rather than geometric in nature are constructed via these theorems.} $$\text{Top}_G\cong\text{Fun}^\times(\text{Fin}_G^\text{op},\text{Top})$$ $$\text{Sp}_G\cong\text{Fun}^\times(\text{Burn}_G,\text{Sp}).$$ Whether a similar result holds for equivariant commutative ring spectra remains an open question. Although we say very little about equivariant homotopy theory in this paper, it is our hope that a more systematic study of $\infty$-categorical Lawvere theories can help us to resolve this question.
\end{example}

We have seen that algebraic theories such as Lawvere theories and PROPs fit very naturally into the framework of categorified commutative algebra, but we have said practically nothing thus far about operads. Operads, though more troublesome to define than Lawvere theories or PROPs, have been central to the development of homotopy theory in recent decades. Since we are using the homotopy-theoretic language of $\infty$-categories and many of our intended applications are to homotopy theory, we would like to do for operads (and $\infty$-operads) what Corollary \ref{CorIntro} does for PROPs and Lawvere theories. The following is a collection of known results translated into our language (see Section 3.4 for details):

\begin{remark}
If $\mathcal{O}$ is an operad, there is an associated symmetric monoidal category $\text{Env}(\mathcal{O})^\odot$, which is a PROP and has the universal property that $$\text{Fun}^\odot(\text{Env}(\mathcal{O})^\odot,\mathcal{C}^\odot)\cong\text{Alg}_{\mathcal{O}}(\mathcal{C}^\odot),$$ for any symmetric monoidal category $\mathcal{C}^\odot$. Moreover, $\text{Env}(\mathcal{O})^\odot\otimes\text{Fin}^\amalg\cong\text{Fin}$. We say that $\text{Env}(\mathcal{O})^\odot$ is \emph{trivial over Fin}.

Similarly, the associated Lawvere theory $\text{Env}(\mathcal{O})^\odot\otimes\text{Fin}^{\text{op},\amalg}$ is trivial over $\text{Burn}^\text{eff}$.
\end{remark}

\begin{conjecture}
The category of operads is equivalent to the full subcategory of Lawvere theories spanned by those which are trivial over $\text{Burn}^\text{eff}$.
\end{conjecture}

In Section 5, we give some evidence for this conjecture, and also describe how to compute the Lawvere theory associated to an operad.

\subsection{Why $\infty$-categories?}
Although all of our results hold for ordinary categories, we use the language of $\infty$-categories throughout. For example, Theorem \ref{IntroThm1} will remain true if we replace every instance of the word `category' by `$\infty$-category' and repeat the rest of the theorem statement verbatim.

We do \emph{not} generalize to $\infty$-categories only out of a desire for needless abstraction, but for three significant reasons:
\begin{enumerate}
\item it considerably simplifies many definitions and (especially) proofs, including that of Theorem \ref{IntroThm1};
\item we can apply our results to situations of homotopy theoretic interest;
\item it is necessary even to formulate Theorem \ref{IntroThm2}.
\end{enumerate}
Point (3) is self-explanatory, but we will elaborate on (1) and (2) shortly. Before we do so, we want to assure the wary reader that we hardly ever have a need to use the technical foundations of $\infty$-categories. Actually, many statements we will need are familiar statements from ordinary category theory which remain true for $\infty$-categories. And in most cases when the results are unfamiliar, we believe they are enlightening rather than technical.

In ordinary commutative algebra, there are two types of perspectives that are available to us. On one hand, we may think of rings as structured sets and manipulate them explicitly. This is typically the easiest way to handle ordinary rings, but we run into problems when we attempt to categorify. There is a tremendous number of axioms which is necessary just to define a commutative semiring category. Even more are needed to make sense of statements like `$\mathcal{C}$ is a module over the semiring category $\mathcal{R}$;' it is possible to write down such definitions, but not at all pleasant to do anything with them. Moreover, the proof (even the statement) of Theorem \ref{IntroThm1} becomes ad hoc and unenlightening.

On the other hand, we may think of rings according to how they act on other rings and on modules. This is the perspective championed by Grothendieck, and it is the perspective suggested by Theorem \ref{IntroThm1}. But this perspective is even more problematic to categorify, since it requires us to make sense of objects like `the 2-category of $\mathcal{R}$-modules, with closed symmetric monoidal structure given by relative tensor product.' We quickly find ourselves pushed up the hierarchy of higher category, to 2-categories and even 3-categories.

Our solution is to work exclusively with $\infty$-categories. The downside of this decision is that the first perspective (the point-set perspective) is very nearly hopeless. But the benefit is that the second perspective (the universal property perspective) is not difficult at all. And unlike with ordinary categories, we are not sucked higher and higher up the hierarchy of $n$-categories. In fact, the $\infty$-category of $\infty$-categories admits the structure of a (large) symmetric monoidal $\infty$-category $\text{Cat}_\infty$, and commutative monoids in $\text{Cat}_\infty$ are just symmetric monoidal $\infty$-categories, and so forth. At no point in this story do we need to utter the dreaded words $(\infty,2)$-category.

In cases where we need to make reference to point-set models (as in Lemma \ref{gpoid}), it is often possible to take homotopy categories and work with the $\infty$-category of 1-categories, where both the point-set and universal property perspectives are available. This is at least an approximation to the $\infty$-categorical world.

As for point (2), some of the most successful applications of algebraic theories (particularly operads) are to homotopy theory, but in this case we typically want to take algebra objects internal to the symmetric monoidal $\infty$-categories of spaces or spectra. Famously, iterated loop spaces and $\mathbb{E}_n$-ring spectra can be described as algebras over little cubes operads \cite{LittleCubes}. In order to make sense of these applications, it will inevitably be necessary to generalize to $\infty$-categories.

In addition, there are other types of algebraic theories which arise in homotopy theory and are not well addressed by the frameworks of operads or Lawvere theories, yet can be understood as cyclic modules over semiring $\infty$-categories. We have already discussed examples from equivariant homotopy theory in Example \ref{IntroEx}. These will not play much of a role in this paper, but they are a primary motivation for the author.

\subsection{Outline}
In Section 2, we set up the background and foundations of the commutative algebra of symmetric monoidal and semiring $\infty$-categories. Very little in this section is original, but nonetheless many of the ideas are likely unfamiliar to anyone but a higher category theory specialist. We are fortunate that we hardly ever need to rely on manipulating fibrations of simplicial sets (which lie at the heart of $\infty$-categories), but are instead able to rely on high-powered \emph{algebraic} tools. As a result, we recommend this section to anyone who wants to see some of the benefits of higher category theory without dwelling on the technical details.

Section 3 is mostly concerned with proving Theorem \ref{IntroThm1} but also contains related results involving the comparison of symmetric monoidal to (co)cartesian monoidal categories -- most notably, a characterization of the `tensoring up' operations $-\otimes\text{Fin}^\amalg$ and $-\otimes\text{Fin}^{\text{op},\amalg}$. And the last part (Section 3.4) covers the examples arising from the study of $\infty$-operads.

In Section 4, we prove Theorem \ref{IntroThm2}, and discuss the extent to which similar statements hold for 1-categories. This provides the comparison between symmetric monoidal $\infty$-categories and spectra.

Finally, in Section 5, we discuss a few conjectures and questions related to $\infty$-operads and computations in categorified commutative algebra.

This paper has changed considerably from an earlier draft, which contained only parts of Sections 2 and 3 and none of Sections 4 or 5. However, the earlier draft contained in addition a more in-depth treatment of Lawvere theories, as providing a bridge between our commutative algebra of symmetric monoidal $\infty$-categories and Lurie's commutative algebra of presentable $\infty$-categories (the two approaches to Principle \ref{Princ} proposed above). It also contained a proof of the last part of Theorem \ref{IntroThm1} (on additive categories). This material is still available on the arXiv (and is complete), but the author hopes to rewrite and expand it in a sequel to this paper \cite{Berman2}, currently in preparation.

\section{Categorified commutative algebra}
In this section, we set up the foundations of the theory of semiring $\infty$-categories and modules over them. There are few new results (with the exception of some in Subsection 2.4), but many of the ideas are not well-known.

In 2.1, we set some notation we will use throughout the rest of the paper and review a few key results from \cite{HA}. Starting in 2.2, we will be extending known (but not always standard) constructions from commutative algebra.

Beginning with the category Set, we may take abelian group objects, and the resulting category Ab is endowed with a closed symmetric monoidal tensor product. Commutative monoid objects in $\text{Ab}^\otimes$ are commutative rings. Everything about this construction is formal, so that Set can be replaced by any presentable $\infty$-category which is closed cartesian monoidal. This is the content of 2.2. When we replace Set by $\text{Cat}_\infty$ (or Cat), we recover commutative semiring $\infty$-categories (or 1-categories) and a theory of commutative algebra over them.

For some commutative rings $R$, being an $R$-module is a \emph{property} of an abelian group, and not extra structure. For example, $\mathbb{Z}/p$-modules correspond to abelian groups which are annihilated by $p$, and $\mathbb{Z}[\frac{1}{p}]$-modules correspond to abelian groups for which multiplication by $p$ is invertible. Such rings are called \emph{solid}, and can be characterized as epimorphisms out of $\mathbb{Z}$, or subsets of $\text{Spec}(\mathbb{Z})$. In 2.3, we generalize this theory, and review known classifications of solid rings and solid ring spectra. This will be useful for us later, as all the main examples of commutative semiring $\infty$-categories we consider in Sections 3 and 4 are solid.

Our most fruitful source of examples of symmetric monoidal $\infty$-categories will arise out of the study of algebraic theories. If we wish to endow an object $X$ of a symmetric monoidal $\infty$-category $\mathcal{C}^\odot$ with some algebraic structure, typically there will be some symmetric monoidal $\infty$-category $\mathcal{T}^\odot$ which encodes all the data of this algebraic structure, and endowing $X$ with the structure corresponds to producing a symmetric monoidal functor $\mathcal{T}^\odot\rightarrow\mathcal{C}^\odot$ which sends a distinguished object of $\mathcal{T}$ to $X$. In order to make sense of the principle `$\mathcal{T}^\odot$ encodes algebraic structure on a \emph{single} object of $\mathcal{C}^\odot$', the distinguished object of $\mathcal{T}^\odot$ should generate all other objects under $\odot$.

That is, algebraic theories will correspond to symmetric monoidal $\infty$-categories which are \emph{cyclic} (generated by one object). In 2.4, we set up a general theory of cyclic modules over a semiring $\infty$-category which will provide the foundation for future study of PROPs, Lawvere theories, and $\infty$-operads (some of which appears in Section 3, and some in the author's upcoming paper \cite{Berman2}).

\begin{remark}
\label{generalize}
We should keep in mind when using ordinary commutative algebra for motivation that we are not only categorifying, but also generalizing in two additional directions:
\begin{enumerate}
\item we are passing to a derived setting (that of $\infty$-categories rather than 1-categories);
\item we are considering semirings instead of rings -- that is, our symmetric monoidal $\infty$-categories do not necessarily have inverses.
\end{enumerate}

For point (1), we should be imagining we are categorifying not commutative rings but rather commutative \emph{ring spectra}. Fortunately, this is not such a problem; the higher algebra and derived geometry of spectra is a vibrant area of current research, and in some respects, ring spectra even behave better than rings (see Remark \ref{SolidRingSp}).

Point (2) is more devastating. As far as the author is aware, the commutative algebra and algebraic geometry of commutative semirings is not very well understood.

Unfortunately, there is no getting around this problem -- practically none of the examples of symmetric monoidal $\infty$-categories that interest us are grouplike (have additive inverses up to equivalence). In fact, we will see in Corollary \ref{CorGpoid} that every commutative \emph{ring} $\infty$-category arises from a commutative ring spectrum.
\end{remark}

\subsection{Background and Notation}
When we say `$\infty$-category', we mean any of the various equivalent definitions, but when in doubt, it is safe to assume we mean the quasicategories of Joyal and Lurie. As much as possible, we follow Lurie's notation from \cite{HTT} and \cite{HA}, and (unless otherwise specified) everything in Section 2.1 is contained in those two sources.

We use Top for the $\infty$-category of Kan complexes (or $\infty$-groupoids), $\text{Sp}$ for spectra, $\text{Cat}_\infty$ for $\infty$-categories, and $\text{SymMon}_\infty$ for symmetric monoidal $\infty$-categories. We also use $\text{CMon}_\infty=\text{CMon}_\infty(\text{Top})$ for $\mathbb{E}_\infty$-spaces (commutative monoids in the $\infty$-category of spaces) and $\text{Ab}_\infty$ for grouplike $\mathbb{E}_\infty$-spaces.

When we write $X\cong Y$, we always mean $X$ and $Y$ are equivalent as objects of some $\infty$-category $\mathcal{C}$. For example, if $X$ and $Y$ are categories (so that $\mathcal{C}=\text{Cat}_1$), we are referring to an equivalence of categories, not an isomorphism. If $\mathcal{C}$ is a 1-category, equivalences and isomorphisms agree, so we may refer to $X\cong Y$ as an isomorphism in that context.

Grouplike $\mathbb{E}_\infty$-spaces are infinite loop spaces, and there is an equivalence of $\infty$-categories $\Omega^\infty\colon\text{Sp}_{\geq 0}\rightarrow\text{Ab}_\infty$ between these and connective spectra. We (somewhat abusively) sometimes refer to the objects of $\text{Ab}_\infty$ themselves as connective spectra.

Without exception, all of the $\infty$-categories just mentioned have natural semiring structures, with two operations $\oplus$ and $\otimes$, and $\otimes$ distributing over $\oplus$ (see 2.2 for details). But when $\oplus$ is cocartesian monoidal, we may write it $\amalg$ instead, and when $\otimes$ is cartesian monoidal, we may write it $\times$ instead. The result is that we will use the notation $\otimes$ in a wide variety of different contexts. To mitigate confusion, we typically use $\odot$ for a symmetric monoidal operation that is not part of a semiring structure, and (as much as possible) indicate which symmetric monoidal structure we are using by a superscript (for example, $\text{Set}^\amalg$ vs. $\text{Set}^\times$). This has the potential to conflict with Lurie's notation in \cite{HA}, where $\mathcal{C}$ and $\mathcal{C}^\otimes$ refer to \emph{different $\infty$-categories}. When there is a possibility of confusion, we will use $[\mathcal{C}^\otimes]$ to refer to the $\infty$-category Lurie calls $\mathcal{C}^\otimes$ (roughly, a category of tuples of objects of $\mathcal{C}$).

If $\mathcal{C}$ is a 1-category, we also use the same notation $\mathcal{C}$ to denote its nerve, regarding $\text{Cat}$ as a full subcategory of $\text{Cat}_\infty$. Thus we think of the 1-categories Set, Fin (finite sets), $\text{Fin}_\ast$ (pointed finite sets), and so forth as $\infty$-categories.

Although constructions involving $\infty$-categories can be cumbersome, most familiar results from ordinary category theory carry over. We will use extensively the following:

\begin{remark}
An $\infty$-category $\mathcal{C}$ with finite products (respectively coproducts) may be endowed with the structure of a symmetric monoidal $\infty$-category, with the symmetric monoidal operation given by products (respectively coproducts). Such a symmetric monoidal $\infty$-category is called \emph{cartesian monoidal} (respectively \emph{cocartesian monoidal}), and we write $\mathcal{C}^\times$ ($\mathcal{C}^\amalg$). A precise definition is given in \cite{HA} 2.4.

Being cartesian (cocartesian) monoidal is a property of a symmetric monoidal $\infty$-category; that is, the forgetful functor $$\text{CartMon}_\infty\rightarrow\text{SymMon}_\infty$$ is a full subcategory inclusion. At the same time, the forgetful functor $$\text{CartMon}_\infty\rightarrow\text{Cat}_\infty$$ is also a subcategory inclusion (equivalent to the subcategory spanned by $\infty$-categories with finite products, and product preserving functors between them).

Given a symmetric monoidal $\infty$-category $\mathcal{C}^\odot$, we can check whether $\mathcal{C}^\odot$ is (co)cartesian monoidal just by checking whether its homotopy category is (co)cartesian monoidal in the ordinary 1-categorical sense (\cite{HA} 2.4.3.19).
\end{remark}

In short, (co)cartesian monoidal $\infty$-categories are far easier to work with than arbitrary symmetric monoidal $\infty$-categories (which are defined to be $\infty$-operads fibered over Comm, satisfying additional properties).

\begin{definition}
A \emph{semiadditive $\infty$-category} is a symmetric monoidal $\infty$-category which is both cocartesian monoidal and cartesian monoidal.
\end{definition}

We are departing from standard terminology in regarding a semiadditive $\infty$-category as a type of symmetric monoidal $\infty$-category rather than a type of $\infty$-category. The latter structure is called preadditive by Gepner, Groth, and Nikolaus (\cite{GGN} Definition 2.1). However, the two definitions are equivalent: an $\infty$-category $\mathcal{C}$ is preadditive if and only if its homotopy category $h\mathcal{C}$ is preadditive/semiadditive (\cite{GGN} Example 2.2), which is true if and only if $h\mathcal{C}$ has cocartesian monoidal and cartesian monoidal structures which agree.

The reader who is not familiar with the higher algebra literature may be especially confused by the discrepancy between the notation $\text{CMon}(\mathcal{C})$ and $\text{CAlg}(\mathcal{C}^\odot)$. By a commutative monoid in an $\infty$-category $\mathcal{C}$ (the former), we mean a product-preserving functor $\text{Fin}\rightarrow\mathcal{C}$; by a commutative algebra in a symmetric monoidal $\infty$-category (the latter), we mean an algebra over the commutative operad. The two notions agree when $\mathcal{C}^\odot$ is cartesian monoidal (so $\odot=\times$).

Following Lurie, we write $\text{Pr}^\text{L}_\infty$ for the $\infty$-category of presentable $\infty$-categories, and left adjoint functors between them. By the adjoint functor theorem, the morphisms of $\text{Pr}^\text{L}_\infty$ can also be regarded as colimit-preserving functors.

There is a closed symmetric monoidal structure on $\text{Pr}^\text{L}$ given by $$\mathcal{C}\otimes\mathcal{D}=\text{Fun}^\text{R}(\mathcal{C}^\text{op},\mathcal{D}),$$ where $\text{Fun}^\text{R}$ denotes right adjoint functors (functors that admit left adjoints).

Finally, commutative algebras in $\text{Pr}^{\text{L},\otimes}$ coincide with closed symmetric monoidal presentable $\infty$-categories (\cite{HA} 4.8). This is also a corollary of the adjoint functor theorem. All of $\text{Top}$, $\text{Sp}$, $\text{CMon}_\infty$, $\text{Ab}_\infty$, $\text{Cat}_\infty$, and $\text{SymMon}_\infty$ are presentable $\infty$-categories, each with a closed symmetric monoidal product (that is, a commutative algebra structure in $\text{Pr}^{\text{L},\otimes}$), which we typically write $\otimes$ instead of $\wedge$ (smash) or $\times$.

\subsection{Semiring and module categories}
We have just asserted that $\text{SymMon}_\infty$ is an example of a presentable $\infty$-category with a closed symmetric monoidal tensor product, but this is certainly not well-known. We believe it was first observed by Gepner, Groth, and Nikolaus \cite{GGN}, as an example of the following lemma:

\begin{lemma}
\label{EasyLemma}
Let $\mathcal{V}^\times$ be a presentable $\infty$-category which is closed cartesian monoidal. Then $\text{CMon}(\mathcal{V})$ has a closed symmetric monoidal structure which is uniquely characterized by the property that the free functor $$\mathcal{V}^\times\rightarrow\text{CMon}(\mathcal{V})^\otimes$$ is symmetric monoidal.

Moreover, if $\mathcal{W}$ is another such $\infty$-category and $L\colon\mathcal{V}\rightleftarrows\mathcal{W}\colon R$ an adjoint pair such that $L$ is product-preserving, there is an induced adjunction $$L_\ast\colon\text{CMon}(\mathcal{V})\rightleftarrows\text{CMon}(\mathcal{W})\colon R_\ast$$ such that:
\begin{enumerate}
\item the adjunction lifts to a symmetric monoidal adjunction (that is, $L_\ast$ lifts to a symmetric monoidal functor and $R_\ast$ to a lax symmetric monoidal functor);
\item $R_\ast$ agrees with $R$ after forgetting commutative monoid structures.
\end{enumerate}
\end{lemma}

Most of the proof is in \cite{GGN}, but the main idea is so striking that we cannot help but include it.

\begin{proof}
Gepner, Groth, and Nikolaus show that $\text{CMon}(\mathcal{V})\cong\text{CMon}_\infty\otimes\mathcal{V}$, where this $\otimes$ is taken in $\text{Pr}^\text{L}_\infty$. Since $\text{CMon}_\infty$ and $\mathcal{V}$ each admit the structure of a commutative algebra object in $\text{Pr}^\text{L}_\infty$ (via closed symmetric monoidal operations $\otimes$ and $\times$, respectively), $\text{CMon}(\mathcal{V})$ inherits such a structure as well, with the given universal property.

Given an adjoint pair $L\colon\mathcal{V}\rightleftarrows\mathcal{W}\colon R$ such that $L$ is product-preserving, $L$ lifts to a functor $\mathcal{V}^\times\rightarrow\mathcal{W}^\times$ in $\text{CAlg}(\text{Pr}^{\text{L},\otimes})$; that is, a symmetric monoidal left-adjoint functor. Tensoring with $\text{CMon}_\infty$, we obtain another symmetric monoidal left-adjoint functor $L_\ast\colon\text{CMon}(\mathcal{V})^\otimes\rightarrow\text{CMon}(\mathcal{W})^\otimes$. By \cite{Haugseng} A.5.11, we have a symmetric monoidal adjunction between $L_\ast$ and $R_\ast$.

We have a commutative diagram of left adjoint functors $$\xymatrix{
\mathcal{V}\ar[r]^L\ar[d]_{\text{Free}} &\mathcal{W}\ar[d]^{\text{Free}} \\
\text{CMon}(\mathcal{V})\ar[r]_{L_\ast} &\text{CMon}(\mathcal{W}).
}$$ Taking right adjoints, we see that $R_\ast$ is compatible with $R$ after forgetting the commutative monoid structures.
\end{proof}

\begin{example}
Taking $\mathcal{V}=\text{Cat}_\infty$, we find that $\text{SymMon}_\infty$ has a symmetric monoidal tensor product characterized by the property that the free functor $\text{Cat}_\infty^\times\rightarrow\text{SymMon}_\infty^\otimes$ is symmetric monoidal.

The same is true for symmetric monoidal 1-categories, taking $\mathcal{V}=\text{Cat}$.
\end{example}

\begin{definition}
A \emph{commutative semiring $\infty$-category} is a commutative algebra object in the symmetric monoidal $\infty$-category $\text{SymMon}_\infty^\otimes$, and we write $$\text{CRingCat}_\infty=\text{CAlg}(\text{SymMon}_\infty^\otimes).$$ For $\mathcal{R}\in\text{CRingCat}_\infty$, we write $\text{Mod}_\mathcal{R}$, $\text{Alg}_\mathcal{R}$, and $\text{CAlg}_\mathcal{R}$ for the $\infty$-categories of $\mathcal{R}$-modules, associative $\mathcal{R}$-algebras, and commutative $\mathcal{R}$-algebras, respectively.
\end{definition}

By `$\mathcal{R}$-modules', we will always mean left/right modules as in \cite{HA} 4.5. Since $\mathcal{R}$ is commutative, the two notions are equivalent.

A commutative semiring $\infty$-category can be interpreted as an $\infty$-category $\mathcal{C}$ with two symmetric monoidal operations $\oplus$ and $\otimes$, satisfying axioms (and even extra data) specifying that $\otimes$ distribute over $\oplus$. We may write $\mathcal{C}^{\oplus,\otimes}$ to emphasize this structure, but in nearly every example it will be unambiguous to write $\mathcal{C}$. Of course, this description is not practical for constructing semiring $\infty$-categories! In particular, the axioms encoding distributivity are infinite in number and not feasible to write down. Instead, most of our examples come from universal constructions.

\begin{example}
Since any presentable $\infty$-category admits all small colimits and any morphism of $\text{Pr}^\text{L}$ preserves all small colimits, there is a forgetful functor $\text{Pr}^\text{L}_\infty\rightarrow\text{SymMon}_\infty$ assigning to a presentable $\infty$-category its cocartesian monoidal structure. Moreover, this functor is lax symmetric monoidal, so any closed symmetric monoidal presentable $\infty$-category $\mathcal{C}^{\otimes}$ is an example of a commutative semiring $\infty$-category $\mathcal{C}^{\amalg,\otimes}$. Here the additive structure is coproduct.

For example, $\text{Set}^{\amalg,\times}$, $\text{Top}^{\amalg,\times}$, $\text{Cat}_\infty^{\amalg,\times}$, $\text{CMon}_\infty^{\vee,\otimes}$, $\text{Ab}_\infty^{\vee,\otimes}$, $\text{SymMon}_\infty^{\times,\otimes}$, and $\text{Sp}^{\vee,\otimes}$ are all commutative semiring $\infty$-categories.
\end{example}

\begin{example}
If $\mathcal{C}^{\oplus,\otimes}$ is a commutative semiring $\infty$-category, and $\mathcal{D}\subseteq\mathcal{C}$ is a full subcategory whose objects are closed (up to equivalence) under $\oplus$ and $\otimes$, then $\mathcal{D}^{\oplus,\otimes}$ is a commutative semiring $\infty$-category as well.

For example, the category Fin of finite sets is a commutative semiring $\infty$-category, as is $\text{Fin}_G$ (finite $G$-sets, for a finite group $G$), as well as the full subcategory of Sp consisting of finite wedges of the sphere spectrum -- this last object is the Burnside $\infty$-category (Example \ref{BurnSub}).
\end{example}

Consider the following adjunctions (the first of which is an equivalence): $$(-)^\text{op}\colon\text{Cat}_\infty\rightleftarrows\text{Cat}_\infty\colon(-)^\text{op}$$ $$i\colon\text{Top}\rightleftarrows\text{Cat}_\infty\colon(-)^\text{iso}$$ $$h\colon\text{Cat}_\infty\rightleftarrows\text{Cat}\colon N.$$ The functor $i$ is the inclusion of $\infty$-groupoids into $\infty$-categories, $N$ is the nerve, and $h$ is the formation of homotopy categories. In both cases, the left adjoint is product-preserving, so Lemma \ref{EasyLemma} applies, producing lax symmetric monoidal functors $(-)^\text{op},i(-)^\text{iso}\colon\text{SymMon}_\infty^\otimes\rightarrow\text{SymMon}_\infty^\otimes$. Taking commutative algebra objects, we learn:

\begin{example}
\label{2.9}
If $\mathcal{C}$ is a commutative semiring $\infty$-category, 
so are $\mathcal{C}^\text{op}$ and $\mathcal{C}^\text{iso}$. The homotopy category $h\mathcal{C}$ is a commutative semiring (1-)category.

For example, $\text{Fin}^\text{iso}$ is a commutative semiring $\infty$-category. In fact, since the free functor $\text{Cat}^\times\rightarrow\text{SymMon}^\otimes$ is symmetric monoidal, and $\text{Fin}^\text{iso}$ is the free symmetric monoidal $\infty$-category on one generator, $\text{Fin}^\text{iso}$ is the unit of $\otimes$.
\end{example}

By analogy with ordinary commutative algebra, we will use $\text{Hom}(-,-)$ to denote the internal Hom in $\text{SymMon}_\infty$, $\text{Hom}_\mathcal{R}(-,-)$ for the inherited internal Hom in $\text{Mod}_\mathcal{R}$, and $\otimes_\mathcal{R}$ for the closed symmetric monoidal relative tensor product in $\text{Mod}_\mathcal{R}$. In particular, $\text{Hom}(-,-)$ will always denote a symmetric monoidal $\infty$-category. For just the mapping space (or mapping $\infty$-groupoid) we may use $\text{Map}_{\text{SymMon}_\infty}(-,-)$ instead.

\begin{remark}
\label{HomUp}
Let $\mathcal{R}\rightarrow\mathcal{A}$ be a functor of semiring $\infty$-categories. The forgetful functor $\text{Mod}_\mathcal{A}\rightarrow\text{Mod}_\mathcal{R}$ has a left adjoint given by $\mathcal{A}\otimes_\mathcal{R} -$ (\cite{HA} 4.5.3), in the sense that it is equivalent to $\mathcal{A}\otimes_\mathcal{R} -$ after restricting along $\text{Mod}_\mathcal{A}\rightarrow\text{Mod}_\mathcal{R}$.

Moreover, the forgetful functor also has a right adjoint $\text{Hom}_\mathcal{R}(\mathcal{A},-)$. To see this, note that $\text{Mod}_\mathcal{A}\rightarrow\text{Mod}_\mathcal{R}$ preserves colimits (\cite{HA} 3.4.4), so by the adjoint functor theorem, it has a right adjoint. The restriction of the right adjoint along $\text{Mod}_\mathcal{A}\rightarrow\text{Mod}_\mathcal{R}$ is itself right adjoint to $\mathcal{A}\otimes_\mathcal{R} -$, so indeed it takes the form $\text{Hom}_\mathcal{R}(\mathcal{A},-)$.
\end{remark}

Note that all of the discussion of this section so far carries over if we begin with the $\infty$-category Cat of categories (or specifically, the nerve of the (2,1)-category Cat). We use the same notation, removing the subscripted $\infty$.

\subsection{Solid rings, ring spectra, and semiring categories}
Throughout Sections 3 and 4, we will give a number of examples of semiring $\infty$-categories $\mathcal{R}$ for which being an $\mathcal{R}$-module is a \emph{property} of a symmetric monoidal $\infty$-category, and not extra structure. There are a number of equivalent ways to formulate this condition.

\begin{definition}
\label{DefSolid}
Let $\mathcal{C}^\odot$ be symmetric monoidal. For $R\in\text{CAlg}(\mathcal{C}^\odot)$, all of the following are equivalent, in which case we call $R$ \emph{solid}.
\begin{enumerate}
\item the forgetful functor $\text{Mod}_R\rightarrow\mathcal{C}$ is fully faithful;
\item the functor $-\odot R\colon\mathcal{C}\rightarrow\mathcal{C}$ is a localization (called a \emph{smashing localization});
\item the multiplication map $R\odot R\rightarrow R$ is an equivalence in $\mathcal{C}$;
\item either of the maps $R\rightarrow R\odot R$ induced by the unit map $1\rightarrow R$ is an equivalence in $\mathcal{C}$;
\item the map $1\rightarrow R$ is an epimorphism in $\text{CAlg}(\mathcal{C}^\odot)$.
\end{enumerate}
\end{definition}

Although we are not aware of a specific work which includes all of these conditions under the name `solid', none of the conditions are new. So the reader who objects that this `definition' requires proof (that all the conditions are equivalent) may consult \cite{GGN} for all but the last condition.

As for (5), it is directly equivalent to (3), because $X\rightarrow Y$ is an epimorphism if and only if (by definition) the codiagonal $Y\amalg_X Y\rightarrow Y$ is an equivalence. The unit 1 is initial in $\text{CAlg}(\mathcal{C}^\odot)$ and so $Y\amalg_1 Y\cong Y\amalg Y$. Moreover, the coproduct $\amalg$ in $\text{CAlg}(\mathcal{C}^\odot)$ is just $\odot$, so the equivalence follows directly.

\begin{remark}
The terminology comes from the classical situation of `solid rings', which were studied and classified by Bousfield and Kan \cite{Core}. The (finitely generated) solid rings are just quotients and localizations of $\mathbb{Z}$, as well as products $\mathbb{Z}[S^{-1}]\times\mathbb{Z}/n$, where each prime divisor of $n$ is in $S$.
\end{remark}

\begin{remark}
\label{georemark}
In commutative rings (and certainly semiring $\infty$-categories), we should not expect epimorphisms to look anything like literal surjections. For example, as in the previous remark, localizations $R\rightarrow R[S^{-1}]$ are epimorphisms. On the other hand, it is true that any surjection is an epimorphism.

Instead, we might think of epimorphisms as having geometric meaning. Provided that $f\colon R\rightarrow A$ is finitely generated, $f$ is an epimorphism if and only if $\text{Spec}(A)\rightarrow\text{Spec}(R)$ is injective in the sense that every fiber is either an isomorphism or empty (\cite{EGA} 17.2.6). So we may think of solid semiring $\infty$-categories as being related to subobjects of the hypothetical geometric `$\text{Spec}(\text{Fin}^\text{iso})$', and therefore telling us something about the geometry of symmetric monoidal $\infty$-categories.
\end{remark}

We would like to port as many techniques from commutative algebra and algebraic geometry as possible into the setting of symmetric monoidal and semiring $\infty$-categories. For example, we might like to classify solid semiring $\infty$-categories. (Although we give a host of examples, a full classification seems out of reach for now.)

But if we want to have any hope of doing this, Remark \ref{generalize} warns us (1) that we should first understand not only solid rings, but solid semirings and solid ring spectra, and (2) that we should expect the problem to arise in generalizing rings to semirings. Indeed, the classification of solid ring spectra is even simpler than that of solid rings. See the next remark for details. Analogously, we will see repeatedly (particularly in Section 2.5) that symmetric monoidal 1-categories can actually be \emph{more} badly behaved than $\infty$-categories.

But solid semirings are very poorly behaved. The author hopes to partially address this problem in an upcoming paper \cite{SemiringBerman}.

\begin{remark}
\label{SolidRingSp}
In the derived setting of ring spectra, epimorphisms have even less in common with surjections. Consider the map of Eilenberg-Maclane ring spectra $\phi\colon HR\rightarrow HA$ induced by a ring map $R\rightarrow A$. For $\phi$ to be an epimorphism, we would need $HA\wedge_{HR} HA\rightarrow HA$ to be an equivalence; that is, not only is $R\rightarrow A$ a ring epimorphism ($A\otimes_R A\rightarrow A$ is an isomorphism), but also $\text{Tor}_\ast^R(A,A)\cong 0$ for all $\ast>0$. This is in general \emph{not} true when $R\rightarrow A$ is surjective.

For example, $H\mathbb{Z}\rightarrow HR$ is an epimorphism of ring spectra only when $R$ is a subring of $\mathbb{Q}$ (localization of $\mathbb{Z}$). While Bousfield and Kan do not prove this, it is a straightforward corollary of their results \cite{Core}, and details can be found in the MathOverflow answer \cite{MathOverflow}.

Note that this is a classification of solid $H\mathbb{Z}$-algebras, not solid ring spectra. But the classification of solid ring spectra is very similar. In particular, a commutative ring spectrum $E$ is solid if and only if it is a Moore spectrum and $\pi_0 E$ is isomorphic to a subring of $\mathbb{Q}$ \cite{SolidSp}.
\end{remark}

\subsection{Cyclic modules}
We now set up the theory of cyclic modules over a semiring $\infty$-category. Unlike in the rest of this section, the classical story of cyclic modules over a ring provides poor motivation here. Classically, cyclic $R$-modules are all quotients of $R$, and are therefore in bijection with ideals of $R$. In fact, they all have ring structures. The parallels break down very quickly. Indeed, if $\mathcal{R}$ is a commutative semiring $\infty$-category, it is usually not even true that cyclic $\mathcal{R}$-modules have semiring structures at all.

Instead, we think of cyclic modules as modeling \emph{algebraic theories}, such as PROPs (Example \ref{PROPExample}), Lawvere theories (\ref{LawvereExample}), or operads (\ref{OperadExample}).

\begin{definition}
Let $\mathcal{R}$ be a commutative semiring $\infty$-category. A \emph{pointed $\mathcal{R}$-module} is an $\mathcal{R}$-module $\mathcal{M}$ along with a choice of distinguished object $X\in\mathcal{R}$; or, equivalently, an $\mathcal{R}$-module $\mathcal{M}$ together with a distinguished map of $\mathcal{R}$-modules $\mathcal{R}\rightarrow\mathcal{M}$. We denote by $$\text{Mod}_{\mathcal{R},\ast}=(\text{Mod}_\mathcal{R})_{\mathcal{R}/}$$ the $\infty$-category thereof.
\end{definition}

\begin{remark}
$\text{Mod}_{\mathcal{R},\ast}$ can be identified with $\mathbb{E}_0$-algebras in $\text{Mod}_\mathcal{R}$ (\cite{HA} 2.1.3.10) and therefore inherits a symmetric monoidal structure $\otimes_\mathcal{R}$ from $\text{Mod}_\mathcal{R}$.
\end{remark}

\begin{definition}
A \emph{cyclic $\mathcal{R}$-module} is a pointed $\mathcal{R}$-module such that the distinguished map $\mathcal{R}\rightarrow\mathcal{M}$ is essentially surjective. We denote by $\text{CycMod}_\mathcal{R}$ the full subcategory of $\text{Mod}_{\mathcal{R},\ast}$ spanned by cyclic modules.
\end{definition}

\begin{example}
\label{PROPExample}
A cyclic $\text{Fin}^\text{iso}$-module (or just cyclic symmetric monoidal $\infty$-category) is a symmetric monoidal $\infty$-category $\mathcal{T}^\odot$ with a distinguished object $X$, such that every object of $\mathcal{T}$ is equivalent to $X^{\odot n}$ for some nonnegative integer $n$. Classically, these are called \emph{PROPs} (product and permutation categories) \cite{PROP}.

Given a PROP $\mathcal{T}^\odot$, a \emph{model} of $\mathcal{T}^\odot$ valued in $\mathcal{C}^\odot$ (often $\mathcal{C}^\odot=\text{Set}^\times$, $\text{Ab}^\otimes$, $\text{Top}^\times$, or $\text{Sp}^\wedge$) is a symmetric monoidal functor $\mathcal{T}^\odot\rightarrow\mathcal{C}^\odot$, and the $\infty$-category thereof is $\text{Mdl}(\mathcal{T}^\odot,\mathcal{C}^\odot)=\text{Hom}(\mathcal{T}^\odot,\mathcal{C}^\odot)$.

When $\mathcal{C}=\text{Top}^\times$, we write just $\text{Mdl}(\mathcal{T}^\odot)=\text{Hom}(\mathcal{T}^\odot,\text{Top}^\times)$.

We can think of the model of $\mathcal{T}^\odot$ as picking out an object of $\mathcal{C}$ (the image of $X$) along with various maps $X^{\odot m}\rightarrow X^{\odot n}$ corresponding to the maps of $\mathcal{T}$. In this way, PROPs model different kinds of algebraic structure (groups, abelian groups, rings, commutative rings, Lie algebras, etc.).
\end{example}

\begin{example}
$\text{Fin}^\amalg$ is the PROP modeling commutative algebras, and similarly $\text{Fin}^{\text{op},\amalg}$ the PROP modeling cocommutative coalgebras.
\end{example}

With this example in mind, we often think of cyclic $\mathcal{R}$-modules as `$\mathcal{R}$-indexed algebraic theories'. This perspective is meaningful even when $\mathcal{R}$ is more complex:

\begin{example}
\label{EqEx}
Fix a finite group $G$, and consider the semiring category of finite $G$-sets with isomorphisms $\text{Fin}_G^\text{iso}$. We may think of cyclic $\text{Fin}_G^\text{iso}$-modules as \emph{equivariant PROPs}.

For example, the subcategory inclusions $\text{Fin}_G^\text{iso}\rightarrow\text{Fin}_G^\text{op}$, $\text{Fin}_G^\text{iso}\rightarrow\text{Burn}_G$ are semiring functors, where $\text{Burn}_G$ is the classical Burnside category of spans of finite $G$-sets. These exhibit $\text{Fin}_G^{\text{op},\amalg}$ and $\text{Burn}_G^\amalg$ are $\text{Fin}_G^\text{iso}$-modules, which are certainly cyclic, and therefore equivariant PROPs.

By Elmendorf's theorem (\cite{Elmendorf} Theorem 1), $\text{Fin}_G^\text{op}$-models in spaces recover the $\infty$-category of genuine equivariant $G$-spaces: $$\text{Hom}(\text{Fin}_G^{\text{op},\amalg},\text{Top}^\times)\cong\text{Top}_G.$$ And by a result of Guillou-May (\cite{GMay2} Theorem 0.1) and Barwick (\cite{BarMack1} Example B.6), $\text{Burn}_G$-models in spectra recover the $\infty$-category of genuine equivariant $G$-spectra: $$\text{Hom}(\text{Burn}_G^\amalg,\text{Sp}^\wedge)\cong\text{Sp}_G.$$

In general, we can think of $\text{Fin}_G^{\text{op},\amalg}$ as the equivariant PROP modeling coefficient system objects, and $\text{Burn}_G$ as the equivariant PROP modeling Mackey functor objects.
\end{example}

\begin{example}
\label{LawvereExample}
If a PROP $\mathcal{T}^\odot$ is cartesian monoidal, it is called a Lawvere theory \cite{Lawvere}. In Example \ref{LawvereExample2}, we will see that Lawvere theories are precisely cyclic $\text{Fin}^\text{op}$-modules.
\end{example}

Although 1-categorical PROPs and Lawvere theories have been studied extensively, higher categorical PROPs and Lawvere theories have only begun to be studied in the past few years. The only sources we are aware of are Cranch's thesis \cite{Cranch} and the appendix of \cite{GGN}.

\begin{example}
\label{OperadExample}
We have already introduced PROPs and Lawvere theories, but \emph{operads} are arguably the most successful model for algebraic theories. We follow Lurie's conventions \cite{HA}, so that an $\infty$-operad $\mathcal{O}$ consists of a type of fibration $[\mathcal{O}^\otimes]\xrightarrow{p}\text{Fin}_\ast$. Remember we are using square brackets to emphasize that $[\mathcal{O}^\otimes]$ is a different $\infty$-category from $\mathcal{O}$. Really, these correspond to what are traditionally called \emph{colored operads}. If we wish to insist that $\mathcal{O}$ be single-colored, we must ask that the underlying $\infty$-category $\mathcal{O}$ have just one object up to equivalence.

By \cite{HA} 2.2.4, to an $\infty$-operad is associated a \emph{symmetric monoidal envelope} $\text{Env}(\mathcal{O})^\odot$. This is a symmetric monoidal $\infty$-category satisfying the universal property $$\text{Hom}(\text{Env}(\mathcal{O})^\odot,\mathcal{C}^\odot)\cong\text{Alg}_{\mathcal{O}}(\mathcal{C}^\odot).$$ If $\mathcal{O}$ is single-colored, then $\text{Env}(\mathcal{O})^\odot$ is a PROP, and its models are exactly $\mathcal{O}$-algebras.

As an $\infty$-category, $\text{Env}(\mathcal{O})$ is the subcategory of $[\mathcal{O}^\otimes]$ spanned by all objects and active morphims between them. That is, it is given by the pullback: $$\xymatrix{
\text{Env}(\mathcal{O})\ar[r]\ar[d] &\text{Fin}\ar[d] \\
[\mathcal{O}^\otimes]\ar[r]_p &\text{Fin}_\ast.
}$$
\end{example}

\begin{example}
\label{EnvExamples}
The commutative ($\mathbb{E}_\infty$) $\infty$-operad has symmetric monoidal envelope Fin (\cite{HA} 2.2.4). The $\mathbb{E}_0$-operad has symmetric monoidal envelope $\text{Fin}^{\text{inj},\amalg}$ (\cite{HA} 2.1.1.19 and 2.2.4).
\end{example}

\begin{remark}
\label{ColoredModule}
In general, if $\mathcal{O}$ is an $\infty$-operad which is not necessarily single-colored, then $\text{Env}(\mathcal{O})^\odot$ is a \emph{colored PROP}. We can make sense of colored algebraic theories using our algebraic language, although they won't play much of a role in this paper.

Specifically, if $\mathcal{R}$ is a commutative semiring $\infty$-category, a \emph{colored $\mathcal{R}$-module} consists of the following data: an $\mathcal{R}$-module $\mathcal{M}$, an $\infty$-category $\mathcal{O}$ of colors, and a fully faithful functor $\mathcal{O}\rightarrow\mathcal{M}$ such that the induced $\mathcal{R}$-module functor $\mathcal{R}[\mathcal{O}]\rightarrow\mathcal{M}$ is essentially surjective. Here, $\mathcal{R}[\mathcal{O}]$ denotes the $\mathcal{R}$-module freely generated by $\mathcal{O}$, which may not be easy to describe though we will not need to do so. However, its objects can be described (up to equivalence) by $\mathcal{R}$-linear combinations of objects of $\mathcal{O}$; that is, formal sums $$\bigoplus_{i=1}^{n}(A_i\otimes X_i)$$ with $A_i\in\mathcal{R}$, $X_i\in\mathcal{O}$.

In other words, an $\mathcal{O}$-colored $\mathcal{R}$-module is an $\mathcal{R}$-module $\mathcal{M}$ and a fully faithful functor $\mathcal{O}\rightarrow\mathcal{M}$ such that every object of $\mathcal{M}$ is equivalent to an $\mathcal{R}$-linear combination of objects in the subcategory $\mathcal{O}$.

For $\mathcal{R}=\text{Fin}^\text{iso}$, respectively $\text{Fin}^\text{op}$, we recover notions of colored PROPs and colored Lawvere theories.
\end{remark}

We end by showing that cyclic modules are well-behaved under relative tensor products.

\begin{lemma}
\label{LemmaES}
Say $\mathcal{R}$ is a commutative semiring $\infty$-category. If $F\colon \mathcal{A}\rightarrow\mathcal{B}$ is an essentially surjective $\mathcal{R}$-module functor and $\mathcal{C}$ is another $\mathcal{R}$-module, then $F_\ast\colon\mathcal{C}\otimes_\mathcal{R}\mathcal{A}\rightarrow\mathcal{C}\otimes_\mathcal{R}\mathcal{B}$ is also essentially surjective.
\end{lemma}

\begin{proof}
Let $\mathcal{S}$ be the full subcategory of $\mathcal{C}\otimes_\mathcal{R}\mathcal{B}$ spanned by the image of $F_\ast$. Then the functor $$j\colon\mathcal{C}\rightarrow\text{Hom}_\mathcal{R}(\mathcal{B},\mathcal{C}\otimes_\mathcal{R}\mathcal{B})$$ (obtained via adjunction from the identity $\mathcal{C}\otimes_\mathcal{R}\mathcal{B}\rightarrow\mathcal{C}\otimes_\mathcal{R}\mathcal{B}$) factors $$\mathcal{C}\rightarrow\text{Hom}_\mathcal{R}(\mathcal{B},\mathcal{S})\subseteq\text{Hom}_\mathcal{R}(\mathcal{B},\mathcal{C}\otimes_\mathcal{R}\mathcal{B}),$$ since the image of $j(X)$ is the same (up to equivalence) as the image of $j(X)\circ F$, but $j(X)\circ F\colon\mathcal{A}\rightarrow\mathcal{C}\otimes_\mathcal{R}\mathcal{B}$ factors through $\mathcal{S}$.

Undoing the tensor-Hom adjunction, the identity on $\mathcal{C}\otimes_\mathcal{R}\mathcal{B}$ factors $$\mathcal{C}\otimes_\mathcal{R}\mathcal{B}\rightarrow\mathcal{S}\subseteq\mathcal{C}\otimes_\mathcal{R}\mathcal{B}$$ (up to equivalence). Therefore, $\mathcal{S}\subseteq\mathcal{C}\otimes_\mathcal{R}\mathcal{B}$ is essentially surjective. Since $\mathcal{S}$ is the image of $F_\ast$, this completes the proof.
\end{proof}

Both of the following propositions are immediate using the lemma.

\begin{proposition}
\label{CyclicTensor}
If $\mathcal{M}$ and $\mathcal{N}$ are two cyclic $\mathcal{R}$-modules, then $\mathcal{M}\otimes_\mathcal{R}\mathcal{N}$ is also cyclic. The structure map from $\mathcal{R}$ is given by $\mathcal{R}\cong\mathcal{R}\otimes_\mathcal{R}\mathcal{R}\rightarrow\mathcal{M}\otimes_\mathcal{R}\mathcal{N}$, the tensor product of the structure maps $\mathcal{R}\rightarrow\mathcal{M},\mathcal{N}$. That is, $\text{CycMod}_\mathcal{R}^{\otimes_\mathcal{R}}$ is symmetric monoidal, and the subcategory inclusion $\text{CycMod}_\mathcal{R}^{\otimes_\mathcal{R}}\subseteq\text{Mod}_{\mathcal{R},\ast}^{\otimes_\mathcal{R}}$ is also symmetric monoidal.
\end{proposition}

\begin{proposition}
For $\mathcal{A}\in\text{CAlg}_\mathcal{R}$, the functor $$\mathcal{A}\otimes_\mathcal{R}-\colon\text{Mod}_\mathcal{R}\rightarrow\text{Mod}_\mathcal{A}$$ restricts to $$\mathcal{A}\otimes_\mathcal{R}-\colon\text{CycMod}_\mathcal{R}\rightarrow\text{CycMod}_\mathcal{A}.$$
\end{proposition}

\begin{corollary}
Higher categorical PROPs form a symmetric monoidal $\infty$-category $\text{PROP}_\infty^\otimes=\text{CycMod}_{\text{Fin}^\text{iso}}^\otimes$ under the ordinary tensor product of symmetric monoidal $\infty$-categories. This tensor product has the following universal property: $$\text{Mdl}(\mathcal{T}\otimes\mathcal{T}^\prime,\mathcal{C}^\odot)\cong\text{Mdl}(\mathcal{T},\text{Mdl}_{\mathcal{T}^\prime}(\mathcal{C}^\odot)).$$
\end{corollary}

We will see in Example \ref{LawvereExample2} that Lawvere theories are also closed under tensor products.

The tensor product of 1-categorical PROPs and Lawvere theories is classical, but we believe we are the first to show that higher PROPs and Lawvere theories form a symmetric monoidal $\infty$-category under $\otimes$.

\section{Cartesian monoidal and semiadditive categories}
Some symmetric monoidal $\infty$-categories $\mathcal{C}^\odot$ may be particularly well-behaved. In a best possible scenario, $\mathcal{C}^\odot$ may be \emph{additive}. That is,
\begin{enumerate}
\item the unit of $\mathcal{C}^\odot$ is a zero object (both initial and terminal);
\item $\odot$ is a categorical biproduct (both a product and a coproduct); and,
\item for every object $X$, there is an automorphism $X\rightarrow X$ corresponding to negation.
\end{enumerate}
Or just some of these conditions may hold: $\mathcal{C}^\odot$ is semi-cartesian monoidal if the unit is terminal or semi-cocartesian monoidal if the unit is initial. If the unit is terminal and $\odot$ is the product (respectively initial and coproduct), $\mathcal{C}^\odot$ is cartesian monoidal (cocartesian monoidal). If (1) and (2) both hold but not necessarily (3), then $\mathcal{C}^\odot$ is semiadditive.

For every such collection of properties, there is a commutative semiring $\infty$-category $\mathcal{R}$ for which being an $\mathcal{R}$-module is equivalent to satisfying that property. Since these are properties and not extra structure, all these semirings will be solid.

In 3.1, we show that Fin and $\text{Fin}^\text{op}$ classify cocartesian and cartesian monoidal $\infty$-categories, and we also show that the effective Burnside 2-category $\text{Burn}^\text{eff}$ classifies semiadditive $\infty$-categories.

In 3.2, we compare symmetric monoidal $\infty$-categories to cocartesian (or, dually, cartesian) monoidal $\infty$-categories via the `tensor up' and `Hom up' operations $\text{Fin}^\amalg\otimes -$ and $\text{Hom}(\text{Fin}^\amalg,-)$.

In 3.3 we discuss the other possible properties from the list above -- with the exception of additive $\infty$-categories themselves, which will be addressed in a sequel work \cite{Berman2}.

An immediate consequence of our results is that PROPs and Lawvere theories can be identified with cyclic modules over $\text{Fin}^\text{iso}$, respectively $\text{Fin}^\text{op}$. Operads also give rise to cyclic modules, but the correspondence is less immediate. In 3.4, we discuss the place of $\infty$-operads in our framework of algebraic theories as cyclic modules. This will be relevant in Section 5, where we present some evidence that there is an equivalence of $\infty$-categories between reduced $\infty$-operads and Lawvere theories which are trivial over $\text{Burn}^\text{eff}$.

\subsection{The semiring categories Fin and $\text{Burn}^\text{eff}$}
We use $\text{CocartMon}_\infty$ (respectively $\text{CartMon}_\infty$) to denote the $\infty$-category of cocartesian (cartesian) monoidal $\infty$-categories, along with those functors that preserve finite coproducts (respectively products).

\begin{theorem}
\label{main1}
$\text{Fin}^{\amalg,\times}$ and $\text{Fin}^{\text{op},\amalg,\times}$ are solid semiring $\infty$-categories. Moreover, a symmetric monoidal $\infty$-category is a Fin-module (respectively $\text{Fin}^\text{op}$-module) if and only if it is cocartesian (respectively cartesian) monoidal. Exactly analogous statements hold for 1-categories.
\end{theorem}

\begin{proof}
We will only prove the claim about Fin-modules for $\infty$-categories; the case of $\text{Fin}^\text{op}$-modules is dual, and the 1-category case is identical, just removing most instances of the symbol $\infty$.

Since Fin is the symmetric monoidal envelope on the commutative operad (Example \ref{EnvExamples}), there is an equivalence $\text{Hom}(\text{Fin}^\amalg,\mathcal{C}^\odot)\xrightarrow{a}\text{CAlg}(\mathcal{C}^\odot)$, with the symmetric monoidal structure on $\text{CAlg}(\mathcal{C}^\odot)$ given by coproducts; that is, the cocartesian monoidal structure (\cite{HA} 3.2.4).

Combining this with \cite{HA} 2.4.3.9, the forgetful functor $\text{CAlg}(\mathcal{C}^\odot)\rightarrow\mathcal{C}^\odot$ is an equivalence if and only if $\mathcal{C}^\odot$ is cocartesian monoidal. Consider the following commutative diagram, with top map induced by inclusion into the second multiplicand $i_2\colon\text{Fin}^\amalg\rightarrow\text{Fin}^\amalg\otimes\text{Fin}^\amalg$: $$\xymatrix{
\text{Map}_{\text{SymMon}_\infty}(\text{Fin}^\amalg,\mathcal{C}^\odot)\ar[d]_a^{\cong} &\text{Map}_{\text{SymMon}_\infty}(\text{Fin}^\amalg\otimes\text{Fin}^\amalg,\mathcal{C}^\odot)\ar[l]_-{f}\ar[d]_{\cong}^a\\
\text{CAlg}(\mathcal{C}^\odot) &\text{CAlg}(\text{CAlg}(\mathcal{C}^\odot)^\amalg)\ar[l]_-{\cong}^-{\text{Forget}}
}$$ Thus $f$ is an equivalence (even of symmetric monoidal $\infty$-categories, but certainly of $\infty$-groupoids). By Yoneda, $i_2\colon\text{Fin}^\amalg\rightarrow\text{Fin}^\amalg\otimes\text{Fin}^\amalg$ is an equivalence of symmetric monoidal $\infty$-categories, and thus Fin is solid.

We now follow an argument as in Section 4.8.2 of \cite{HA} to show that a symmetric monoidal $\infty$-category admits the structure of a Fin-module if and only if it is cocartesian monoidal. This will suffice, since $\text{Mod}_\text{Fin}\rightarrow\text{SymMon}_\infty$ is fully faithful (because Fin is solid).

Suppose that $\mathcal{C}^\amalg$ is a cocartesian monoidal $\infty$-category. Then $$\mathcal{C}\cong\text{CAlg}(\mathcal{C}^\amalg)\cong\text{Hom}(\text{Fin}^\amalg,\mathcal{C}^\amalg),$$ and the latter lifts to a Fin-module structure, as desired (see the description of the right adjoint in Remark \ref{HomUp}).

Conversely, suppose $\mathcal{C}^\odot$ is a Fin-module. Then $$\mathcal{C}^\odot\cong\text{Hom}_\text{Fin}(\text{Fin}^\amalg,\mathcal{C}^\odot)\cong\text{Hom}(\text{Fin}^\amalg,\mathcal{C}^\odot)\cong\text{CAlg}(\mathcal{C}^\odot)^\amalg,$$ which is cocartesian monoidal. (The second equivalence in the chain is because Fin is solid.)
\end{proof}

\begin{example}
\label{LawvereExample2}
In Section 2.4, we asserted that cyclic modules over semiring $\infty$-categories are models for algebraic theories, and we saw that $$\text{PROP}_\infty^\otimes\cong\text{CycMod}_{\text{Fin}^\text{iso}}^\otimes.$$

Cyclic $\text{Fin}^\text{op}$-modules are PROPs which are cartesian monoidal. Classically, these are known as \emph{Lawvere theories}, another common model for algebraic theories. For the symmetric monoidal $\infty$-category thereof, we write $$\text{Lawvere}^\otimes_\infty=\text{CycMod}_{\text{Fin}^\text{op}}^\otimes;$$ the tensor product comes from Proposition \ref{CyclicTensor}.
\end{example}

If $\text{Fin}$ classifies cocartesian monoidal $\infty$-categories, and $\text{Fin}^\text{op}$ classifies cartesian monoidal $\infty$-categories, then $\text{Fin}^\amalg\otimes\text{Fin}^{\text{op},\amalg}$ should classify those symmetric monoidal $\infty$-categories which are both cartesian and cocartesian -- that is, semiadditive $\infty$-categories. A natural problem then is to compute $\text{Fin}\otimes\text{Fin}^\text{op}$. It turns out that it is the 2-category of spans of finite sets, also known as the effective Burnside 2-category.

\begin{definition}
The \emph{effective Burnside 2-category} $\text{Burn}^\text{eff}$ (which we called $\text{Span}(\text{Fin})$ in the introduction) has finite sets for objects, and a morphism from $S$ to $T$ consists of another finite set $X$ and a span $$\xymatrix{
&X\ar[ld]\ar[rd] &\\
S &&T.
}$$ Composition is via pullbacks, and 2-morphisms (all of which are invertible) are bijections of the top object in the span. Via the nerve, we think of $\text{Burn}^\text{eff}$ as an $\infty$-category.
\end{definition}

There is also a more general effective Burnside (or span) construction. The construction for categories is classical; for $\infty$-categories it is due to Barwick \cite{BarMack1}, so we use his language (`effective Burnside construction'). In particular, $\text{Burn}^\text{eff}$ is just the effective Burnside construction applied to Fin. For now, we will not need to understand the details of this construction in any more generality.

\begin{corollary}
\label{main2}
As symmetric monoidal $\infty$-categories, $$\text{Burn}^{\text{eff},\amalg}\cong\text{Fin}^{\text{op},\amalg}\otimes\text{Fin}^\amalg.$$ In particular, $\text{Burn}^{\text{eff},\amalg,\times}$ is a commutative semiring $\infty$-category, with $\times$ distributing over $\amalg$, and it is solid. Moreover, a symmetric monoidal $\infty$-category admits the structure of a $\text{Burn}^\text{eff}$-module if and only if it is semiadditive. Analogous statements hold for 1-categories, with $\text{Burn}^\text{eff}$ replaced by its homotopy category $\text{hBurn}^\text{eff}$.
\end{corollary}

\begin{proof}
We begin with the $\infty$-category case.

We first prove that $\text{Fin}\otimes\text{Fin}^\text{op}$-modules are just semiadditive $\infty$-categories. Indeed, since Fin and $\text{Fin}^\text{op}$ are solid, they are idempotent by Definition \ref{DefSolid}. Thus $\text{Fin}\otimes\text{Fin}^\text{op}$ is also solid, so that $\text{Mod}_{\text{Fin}\otimes\text{Fin}^\text{op}}$ is the full subcategory of $\text{SymMon}_\infty$ given by the intersection of the full subcategories $\text{Mod}_\text{Fin}\cong\text{CocartMon}_\infty$ and $\text{Mod}_{\text{Fin}^\text{op}}\cong\text{CartMon}_\infty$. A semiadditive $\infty$-category is just a cartesian monoidal $\infty$-category which is also cocartesian monoidal, so the claim is certainly true of $\text{Fin}\otimes\text{Fin}^\text{op}$.

Now it suffices to show that $\text{Burn}^{\text{eff},\amalg}\cong\text{Fin}^\amalg\otimes\text{Fin}^{\text{op},\amalg}$ as symmetric monoidal $\infty$-categories. It will follow that $\text{Burn}^\text{eff}$ has all the properties described. Consider the free functors (that is, left adjoints to the forgetful functors) $$\xymatrix{
\text{Cat}_\infty\ar[rr]^-{\text{Fin}^\text{iso}[-]} &&\text{SymMon}_\infty\ar[rr]^-{-\otimes\text{Fin}\otimes\text{Fin}^\text{op}} &&\text{SemiaddCat}_\infty
}$$ The free semiadditive $\infty$-category on the trivial one-object $\infty$-category $\ast$ is $\text{Fin}\otimes\text{Fin}^\text{op}$. 

On the other hand, Glasman (\cite{Glasman} Theorem A.1) has already done our hard work for us by showing that $\text{Burn}^\text{eff}$ is the free semiadditive $\infty$-category on $\ast$. That completes the proof for $\infty$-categories.

The only difference for the 1-category case is that we need to know $\text{hBurn}^\text{eff}$ is the free semiadditive category on one generator. This can be proven directly, or alternatively follows from the $\infty$-categorical statement since the homotopy category functor $\text{h}\colon\text{Cat}_\infty\rightarrow\text{Cat}$ is itself free; that is, the following diagram commutes since the diagram of left adjoints (made up of horizontal forgetful functors and vertical nerves) also commutes: $$\xymatrix{
\text{Cat}_\infty\ar[rr]^-{\text{Free}}\ar[d]_{\text{h}} &&\text{SemiaddCat}_\infty\ar[d]_{\text{h}} \\
\text{Cat}\ar[rr]_-{\text{Free}} &&\text{SemiaddCat}.
}$$
\end{proof}

\begin{remark}
\label{PK}
Since $\text{CocartMon}_\infty\cong\text{Mod}_\text{Fin}$, $\text{CocartMon}_\infty$ inherits a symmetric monoidal product $\otimes_\text{Fin}$. Moreover, since Fin is solid, the tensor product or Hom of two cocartesian monoidal $\infty$-categories agrees with their tensor product or Hom as symmetric monoidal $\infty$-categories: $\mathcal{C}\otimes_\text{Fin}\mathcal{D}\cong\mathcal{C}\otimes\mathcal{D}$ and $\text{Hom}_\text{Fin}(\mathcal{C},\mathcal{D})\cong\text{Hom}(\mathcal{C},\mathcal{D})$.

Lurie has constructed a tensor product for $\infty$-categories closed under certain colimits \cite{HA} 4.8.1. By considering $\infty$-categories closed under finite coproducts, his construction produces a closed symmetric monoidal structure for cocartesian monoidal $\infty$-categories. This agrees with ours, because for each $\mathcal{C}$ admitting finite coproducts, the functors $-\otimes^\text{Lurie}\mathcal{C}$ and $-\otimes\mathcal{C}^\amalg$ have the same right adjoint $\text{Fun}^\amalg(\mathcal{C},-)\cong\text{Hom}(\mathcal{C}^\amalg,-)$ (\cite{HA} 4.8.1.6).
\end{remark}

\subsection{Free and cofree (co)cartesian monoidal categories}
If we wish to understand the relationship between $\text{SymMon}_\infty$ and $\text{Mod}_\mathcal{R}$, where $\mathcal{R}$ is one of the semiring $\infty$-categories just discussed, then we might study the cofree construction $\text{Hom}(\mathcal{R},-)$ and free construction $\mathcal{R}\otimes -$ as in Remark \ref{HomUp}.

\begin{remark}
\label{cofree}
In the cases of $\mathcal{R}=\text{Fin}$, $\text{Fin}^\text{op}$, or $\text{Burn}^\text{eff}$, we already know the effect of the cofree construction: the commutative algebra construction $$\text{CAlg}\cong\text{Hom}(\text{Fin}^\amalg,-)\colon\text{SymMon}_\infty\rightarrow\text{CocartMon}_\infty$$ and the cocommutative coalgebra construction $$\text{CocCoalg}\cong\text{Hom}(\text{Fin}^{\text{op},\amalg},-)\colon\text{SymMon}_\infty\rightarrow\text{CartMon}_\infty$$ are the respective cofree functors from symmetric monoidal $\infty$-categories to (co)cartesian monoidal $\infty$-categories.

The commutative-cocommutative bialgebra construction is equivalent to $$\text{Hom}(\text{Burn}^{\text{eff},\amalg},-)\cong\text{CAlg}(\text{CCoalg}(-)),$$ the cofree functor from symmetric monoidal $\infty$-categories to semiadditive $\infty$-categories.

For us, these results are direct consequences of Theorem \ref{main1} with \ref{HomUp}.
\end{remark}

The free construction is more complicated, and not always easy to compute, but we do have some results.

First, recall that if $\mathcal{C}^\times,\mathcal{D}^\times$ are two cartesian monoidal categories, then $\text{Hom}(\mathcal{C}^\times,\mathcal{D}^\times)$ is the full subcategory of $\text{Fun}(\mathcal{C},\mathcal{D})$ spanned by product-preserving functors. The following is a sort of cartesian monoidal Yoneda lemma.

\begin{lemma}
\label{SMYoneda}
If $\mathcal{C}^\times$ is a cartesian monoidal $\infty$-category, the Yoneda embedding $\mathcal{C}^\text{op}\rightarrow\text{Fun}(\mathcal{C},\text{Top})$ factors through $\text{Hom}(\mathcal{C}^\times,\text{Top}^\times)$, so that there is a full subcategory inclusion $$\mathcal{C}^\text{op}\subseteq\text{Hom}(\mathcal{C}^\times,\text{Top}^\times).$$ Dually, if $\mathcal{C}^\amalg$ is cocartesian monoidal, there is a full subcategory inclusion $$\mathcal{C}\subseteq\text{Hom}(\mathcal{C}^{\text{op},\amalg},\text{Top}^\times).$$ These embeddings are themselves symmetric monoidal functors if and only if $\mathcal{C}$ is semiadditive.
\end{lemma}

\begin{proof}
By \cite{HTT} 5.1.3.2, corepresentable functors $\mathcal{C}\rightarrow\text{Top}$ preserve limits, and in particular finite products. Thus if $\mathcal{C}^\times$ is Cartesian monoidal, then the Yoneda embedding factors through $\mathcal{C}^\text{op}\rightarrow\text{Hom}(\mathcal{C}^\times,\text{Top}^\times)$. The proof is similar for $\mathcal{C}^\amalg$ cocartesian monoidal.

For the last sentence: If (for example) $\mathcal{C}^\times$ is cartesian monoidal, then $$\mathcal{C}^{\text{op},\times}\rightarrow\text{Hom}(\mathcal{C}^\times,\text{Top}^\times)^\times$$ is symmetric monoidal if and only if (applying tensor-Hom adjunctions), $\mathcal{C}\rightarrow\text{Fun}(\mathcal{C}^\text{op},\text{Top})$ factors through $\text{Hom}(\mathcal{C}^{\text{op},\times},\text{Top}^\times)$. That is, every corepresentable functor $\mathcal{C}(-,X)$ sends finite products in $\mathcal{C}$ to products in $\text{Top}$. This is true if and only if $$\mathcal{C}(A\times B,-)\cong\mathcal{C}(A,-)\times\mathcal{C}(B,-)\cong\mathcal{C}(A\amalg B,-),$$ which is to say if and only if $\mathcal{C}$ is semiadditive.
\end{proof}

\begin{example}
If $\mathcal{T}^\times$ is any Lawvere theory, $\mathcal{T}$ is a full subcategory of $$\text{Mdl}(\mathcal{T}^\times)^\text{op}\cong\text{Hom}(\mathcal{T}^\times,\text{Top}^\times)^\text{op}.$$
\end{example}

\begin{remark}
If $\mathcal{C}^\odot$ is symmetric monoidal, then there is a forgetful functor $\text{Fgt}_\mathcal{C}\colon\text{Hom}(\mathcal{C}^{\text{op},\odot},\text{Top}^\times)\rightarrow\text{Fun}(\mathcal{C}^\text{op},\text{Top})$, but in general we should not expect $\text{Fgt}_\mathcal{C}$ to be fully faithful (if $\mathcal{C}^\odot$ is not cocartesian monoidal).
\end{remark}

\begin{theorem}
\label{TensorFin}
If $\mathcal{C}^\odot$ is a symmetric monoidal $\infty$-category, then $\text{Fgt}_\mathcal{C}$ (as in the above remark) has a left adjoint $\text{Free}_\mathcal{C}$. If $X\in\mathcal{C}$, let $\underline{X}$ denote the representable functor $\text{Map}_\mathcal{C}(-,X)$.

Then the free cocartesian monoidal $\infty$-category $\mathcal{C}^\odot\otimes\text{Fin}^\amalg$ generated by $\mathcal{C}^\odot$ is equivalent to the full subcategory of $\text{Hom}(\mathcal{C}^{\text{op},\odot},\text{Top}^\times)$ spanned by the objects $\text{Free}_\mathcal{C}(\underline{X})$, $X\in\mathcal{C}$, with its induced cocartesian monoidal structure.
\end{theorem}

\begin{remark}
Dually (by considering $\mathcal{C}^\text{op}$ instead of $\mathcal{C}$), $\mathcal{C}^\odot\otimes\text{Fin}^{\text{op},\amalg}$ is the opposite of the corresponding full subcategory of $\text{Hom}(\mathcal{C}^\odot,\text{Top}^\times)$.
\end{remark}

\begin{proof}
Denote $\mathcal{P}(\mathcal{C})=\text{Fun}(\mathcal{C}^\text{op},\text{Top})$ and $\mathcal{P}_\odot(\mathcal{C}^\odot)=\text{Hom}(\mathcal{C}^{\text{op},\odot},\text{Top}^\times)$ throughout the proof. We may write $\mathcal{P}_\amalg(\mathcal{C}^\amalg)$ instead of $\mathcal{P}_\odot(\mathcal{C}^\amalg)$ for cocartesian monoidal $\mathcal{C}^\amalg$.

First, assume $\mathcal{C}^\amalg$ cocartesian monoidal, so that $\mathcal{C}^\amalg\otimes\text{Fin}^\amalg\cong\mathcal{C}^\amalg$. By \cite{HTT} 5.3.6.10, the forgetful functor $\text{Pr}^\text{L}_\infty\rightarrow\text{CocartMon}_\infty$ has a left adjoint given by $\mathcal{P}_\amalg(-)$, and if $F\colon\mathcal{C}^\amalg\rightarrow\mathcal{D}^\amalg$ is a coproduct-preserving functor, $\mathcal{P}_\amalg(F)$ is left adjoint to restriction $F^\ast\colon\text{Hom}(\mathcal{D}^{\text{op},\amalg},\text{Top}^\times)\rightarrow\text{Hom}(\mathcal{C}^{\text{op},\amalg},\text{Top}^\times)$.

As in \cite{HTT} 5.3.6 (see also Remark \ref{PK}), let $\mathcal{P}_\emptyset^\amalg$ denote the left adjoint to the forgetful functor $\text{CocartMon}_\infty\rightarrow\text{Cat}_\infty$. (Elsewhere in this paper, we have called $\mathcal{P}_\emptyset^\amalg(\mathcal{C})=\text{Fin}[\mathcal{C}]$.) And let
$F\colon\mathcal{P}_\emptyset^\amalg(\mathcal{C})\rightarrow\mathcal{C}^\amalg$, $G\colon\mathcal{C}\rightarrow\mathcal{P}_\emptyset^\amalg(\mathcal{C})$ be the two functors arising from the unit and counit of the adjunction $(\mathcal{P}_\emptyset^\amalg,\text{Forget})$. Precompositions with $F$ and $G$ induce $$\mathcal{P}_\amalg(\mathcal{C}^\amalg)\xrightarrow{F^\ast}\mathcal{P}_\amalg(\mathcal{P}_\emptyset^\amalg(\mathcal{C}))\xrightarrow{G^\ast}\mathcal{P}(\mathcal{C}).$$ The composite functor is $\text{Fgt}_\mathcal{C}$ (because $F\circ G\colon\mathcal{C}\rightarrow\mathcal{C}$ is the identity), $G^\ast$ is an equivalence by \cite{HTT} 5.3.6.10, and $F^\ast$ is right adjoint to $\mathcal{P}_\amalg(F)$. Therefore $\text{Fgt}_\mathcal{C}$ has a left adjoint $\text{Free}_\mathcal{C}$.

Since $\text{Free}_\mathcal{C}$ has fully faithful right adjoint, it is a localization functor. In particular, the restriction $$\text{Free}_\mathcal{C}\colon\mathcal{P}_\amalg(\mathcal{C}^\amalg)\subseteq\mathcal{P}(\mathcal{C})\rightarrow\mathcal{P}_\amalg(\mathcal{C}^\amalg)$$ is naturally equivalent to the identity (\cite{HTT} 5.2.7.4). The theorem's claim follows from $\text{Free}_\mathcal{C}(\underline{X})\cong\underline{X}$ and Lemma \ref{SMYoneda}.

This completes the proof if $\mathcal{C}^\amalg$ is cocartesian monoidal, but now take $\mathcal{C}^\odot$ arbitrary. The symmetric monoidal functor $i\colon\text{Fin}^\text{iso}\rightarrow\text{Fin}$ (which realizes Fin as cyclic) induces after tensoring with $\mathcal{C}$ a symmetric monoidal functor $$F\colon\mathcal{C}^\odot\rightarrow\mathcal{C}^\odot\otimes\text{Fin}^\amalg$$ which is essentially surjective by Lemma \ref{LemmaES}. Moreover, \\$i^\ast\colon\text{Hom}(\text{Fin}^{\text{op},\amalg},\text{Top}^\times)\rightarrow\text{Hom}(\text{Fin}^{\text{iso},\amalg},\text{Top}^\times)$ is an equivalence, so $$F^\ast\colon\text{Hom}(\mathcal{C}^{\text{op},\odot}\otimes\text{Fin}^{\text{op},\amalg},\text{Top}^\times)\rightarrow\text{Hom}(\mathcal{C}^{\text{op},\odot},\text{Top}^\times)$$ is also an equivalence (from the tensor-Hom adjunction).

Call $\mathcal{D}^\odot=\mathcal{C}^\odot\otimes\text{Fin}^\amalg$, and consider the commutative square\footnote{For the sake of this proof, we only need to know that the squares commute in the weakest sense: the two composite functors around the square should send any single object to equivalent objects.} $$\xymatrix{
\mathcal{P}_\odot(\mathcal{C}^\odot)\ar[r]^-{\text{Fgt}_\mathcal{C}} &\mathcal{P}(\mathcal{C}) \\
\mathcal{P}_\odot(\mathcal{D}^\odot)\ar[r]_-{\text{Fgt}_\mathcal{D}}^-{\subseteq}\ar[u]_{\cong}^{F^\ast} &\mathcal{P}(\mathcal{D})\ar[u]_{F^\ast}.
}$$ We have already established that $\text{Fgt}_\mathcal{D}$ has a left adjoint (because $\mathcal{D}^\odot$ is cocartesian monoidal). In addition, the left vertical functor is an equivalence, and the right vertical functor has left adjoint $\mathcal{P}(F)$ by \cite{HTT} 5.1.5.6. Thus $\text{Fgt}_\mathcal{C}\cong F^\ast\circ\text{Fgt}_\mathcal{D}\circ(F^\ast)^{-1}$ has a left adjoint $\text{Free}_\mathcal{C}$. Taking left adjoints, we have another commutative square (on the left), and the right square also commutes since, by the Yoneda lemma, $F^\ast(\underline{X})\cong\underline{F(X)}$: $$\xymatrix{
\mathcal{P}_\odot(\mathcal{C}^\odot)\ar[d]^{\cong}_{(F^\ast)^{-1}} &\mathcal{P}(\mathcal{C})\ar[l]_-{\text{Free}_\mathcal{C}}\ar[d]^{F_\ast} &\mathcal{C}\ar[l]_-{\subseteq}\ar[d]^F \\
\mathcal{P}_\odot(\mathcal{D}^\odot) &\mathcal{P}(\mathcal{D})\ar[l]^-{\text{Free}_\mathcal{D}} &\mathcal{D}\ar[l]^-{\subseteq}.
}$$ The bottom composite $\mathcal{D}\rightarrow\mathcal{P}_\odot(\mathcal{D}^\odot)$ is fully faithful because $\mathcal{D}^\odot$ is cocartesian monoidal (Lemma \ref{SMYoneda}) and $\text{Free}_\mathcal{D}(\underline{X})\cong\underline{X}$ (first part of the proof).

Since moreover $\mathcal{C}\rightarrow\mathcal{D}$ is essentially surjective, $\mathcal{D}=\mathcal{C}\otimes\text{Fin}$ is the full subcategory of $\mathcal{P}_\odot(\mathcal{C}^\odot)$ spanned by the image of $\mathcal{C}$. Following the top horizontal functors, the image of $X\in\mathcal{C}$ in $\mathcal{P}_\odot(\mathcal{C}^\odot)$ is $\text{Free}_\mathcal{C}(\underline{X})$, as desired.
\end{proof}

Thus computation of the free functor $\text{SymMon}_\infty\rightarrow\text{CocartMon}_\infty$ reduces to computation of another free functor $\text{Fun}(\mathcal{C}^\text{op},\text{Top})\rightarrow\text{Hom}(\mathcal{C}^{\text{op},\odot},\text{Top}^\times)$, which may very well still be a difficult problem.

However, if $\mathcal{C}^\odot$ is cyclic (a PROP), then we can say more.

\begin{example}
The Lawvere theory associated to a PROP $\mathcal{T}^\odot$ is $\mathcal{T}^\odot\otimes\text{Fin}^{\text{op},\amalg}$ (which is a Lawvere theory by Proposition \ref{CyclicTensor}). If $\mathcal{C}^\times$ is cartesian monoidal (for example, Set or Top), then $$\text{Mdl}(\mathcal{T}^\odot\otimes\text{Fin}^{\text{op},\amalg},\mathcal{C}^\times)\cong\text{Mdl}(\mathcal{T}^\odot,\mathcal{C}^\times).$$
\end{example}

\begin{corollary}
\label{TensorFinPROP}
The Lawvere theory $\mathcal{T}^\odot\otimes\text{Fin}^{\text{op},\amalg}$ associated to a PROP $\mathcal{T}^\odot$ is equivalent to the full subcategory of $\text{Mdl}(\mathcal{T}^\odot)^\text{op}$ spanned by the models which are freely generated by finite sets (finitely generated free models).
\end{corollary}

\begin{proof}
Let $X$ be the distinguished (generating) object of $\mathcal{T}$. We have a forgetful functor $$R\colon\text{Hom}(\mathcal{T}^{\text{op},\odot},\text{Top}^\times)\rightarrow\text{Top},$$ which evaluates at $X$. We want to show that the left adjoint $L$ sends the finite discrete space on $n$ points to $\text{Free}_\mathcal{T}(\underline{X^{\odot n}})$, up to equivalence.

But indeed mapping into $F\in\text{Hom}(\mathcal{T}^{\text{op},\odot},\text{Top}^\times)$ from $L(n)$ or $\text{Free}_\mathcal{T}(\underline{X^{\odot n}})$ has the same universal property, $$\text{Map}(L(n),F)\cong F(X)^{\times n}\cong F(X^{\odot n})\cong\text{Map}(\text{Free}_\mathcal{T}(\underline{X^{\odot n}}),F).$$
\end{proof}

\begin{example}
Since $\text{Burn}^\text{eff}$ is the Lawvere theory for $\mathbb{E}_\infty$-spaces, it is the full subcategory of $\text{CMon}_\infty$ spanned by those $\mathbb{E}_\infty$-spaces which are freely generated by finite sets.

By the 1-categorical analogues of these results (which are also true with the same proofs), the 1-category $\text{hBurn}^\text{eff}$ is the full subcategory of commutative monoids spanned by $\mathbb{N}^k$ as $k\geq 0$ varies.
\end{example}

\begin{example}
\label{BurnSub}
Let Burn be the Lawvere theory for $\text{Ab}_\infty$ (which exists by the appendix of \cite{GGN}). Burn is the full subcategory of Sp spanned by wedge powers of the sphere spectrum, because $\text{Burn}\subseteq\text{Hom}(\text{Burn}^{\text{op},\vee},\text{Top}^\times)\cong\text{Ab}_\infty\subseteq\text{Sp}$.

Burn can also be obtained from $\text{Burn}^\text{eff}$ by group-completing the Hom-spaces (which are already $\mathbb{E}_\infty$-spaces since $\text{Burn}^\text{eff}$ is semiadditive). This construction is known to many people, but we are not aware of a reference for $\infty$-categories. It will appear in detail in the sequel to this paper \cite{Berman2}.

Similarly, hBurn is the full subcategory of abelian groups spanned by the finitely generated free abelian groups.
\end{example}

\begin{example}
\label{OperadCompute}
If $\mathcal{O}$ is a (single-colored) $\infty$-operad, the associated Lawvere theory is $\text{Env}(\mathcal{O})^\odot\otimes\text{Fin}^{\text{op},\amalg}$.

By Corollary \ref{TensorFinPROP}, this Lawvere theory is the (opposite of the) full subcategory of $\text{Alg}_\mathcal{O}(\text{Top}^\times)$ spanned by the free $\mathcal{O}^\otimes$-algebras on finite sets. But there is a well-known calculation of free algebras over an operad; the free $\mathcal{O}^\otimes$-algebra on the finite set $n$ has underlying space $$\coprod_{T\in\text{Fin}}(\mathcal{O}(T)\times_{\Sigma_T}n^T).$$ (Here we are using classical language; see \cite{HA} 3.1 for a precise treatment of free algebras over $\infty$-operads.) In particular, we can describe all the mapping spaces in $\text{Env}(\mathcal{O})^\odot\otimes\text{Fin}^{\text{op},\amalg}$: $$\text{Map}(X^{\odot m},X^{\odot n})\cong (\coprod_{T\in\text{Fin}}\mathcal{O}(T)\times_{\Sigma_T}m^T)^n.$$
\end{example}

\subsection{The semiring categories $\text{Fin}^\text{inj}$ and $\text{Fin}_\ast$}
Recall that we are using $\text{Fin}^\text{inj}$ for the category of finite sets, considering only the injections between them.

\begin{proposition}
\label{PropInj}
$\text{Fin}^\text{inj}$ and its opposite are solid semiring $\infty$-categories. A symmetric monoidal $\infty$-category $\mathcal{C}^\odot$ is a $\text{Fin}^\text{inj}$-module (respectively $\text{Fin}^{\text{inj},\text{op}}$-module) if and only if the unit of $\mathcal{C}^\odot$ is initial (respectively terminal). We say that $\mathcal{C}^\odot$ is semi-cocartesian monoidal (respectively semi-cartesian monoidal).
\end{proposition}

\begin{proof}
The proof is exactly as in Theorem \ref{main1}, noting that $\text{Fin}^{\text{inj},\amalg}$ is the symmetric envelope of the $\mathbb{E}_0$-operad (Example \ref{EnvExamples}), and $\text{Alg}_{\mathbb{E}_0}(\mathcal{C}^\odot)\cong\mathcal{C}_{1/}$ (\cite{HA} 2.1.3.10).
\end{proof}

\begin{example}
\label{LFinj}
There is an equivalence of symmetric monoidal $\infty$-categories $\text{Fin}^\amalg\otimes\text{Fin}^{\text{inj},\text{op},\amalg}\cong\text{Fin}_\ast^\amalg$. This follows from Corollary \ref{TensorFinPROP} along with the observation that $$\text{Hom}(\text{Fin}^{\text{inj},\amalg},\text{Top}^\times)\cong\text{Alg}_{\mathbb{E}_0}(\text{Top}^\times)=\text{Top}_\ast.$$
\end{example}

The following proposition follows immediately from Example \ref{LFinj}.

\begin{proposition}
$\text{Fin}_\ast$ and its opposite are solid semiring $\infty$-categories. A symmetric monoidal $\infty$-category $\mathcal{C}$ is a $\text{Fin}_\ast$-module (respectively $\text{Fin}_\ast^\text{op}$-module) if and only if $\mathcal{C}$ is cocartesian (respectively cartesian) monoidal and has a zero object (an object which is both initial and terminal).
\end{proposition}

\begin{remark}
\label{subBurn}
For a given morphism $X\leftarrow T\rightarrow Y$ in $\text{Burn}^\text{eff}$, we will say that the morphism is, for example, (injective,arbitrary) if $T\rightarrow X$ is injective and $T\rightarrow Y$ is arbitrary. By considering functions which are injective, bijective, or arbitrary, we obtain in this way 9 types of structured morphisms in $\text{Burn}^\text{eff}$. For each type of structured morphism, we may consider the subcategory of $\text{Burn}^\text{eff}$ spanned by all objects, but only the morphisms of that type. In this way, we recover each of the 9 semiring $\infty$-categories we have considered so far:
\begin{itemize}
\item $\text{Fin}^\text{inj}$ -- (bijective,injective) morphisms;
\item $\text{Fin}^{\text{inj},\text{op}}$ -- (injective,bijective) morphisms;
\item $\text{Fin}$ -- (bijective,arbitrary) morphisms;
\item $\text{Fin}^\text{op}$ -- (arbitrary,bijective) morphisms;
\item $\text{Fin}_\ast\cong\text{Fin}\otimes\text{Fin}^{\text{inj},\text{op}}$ -- (injective,arbitrary) morphisms;
\item $\text{Fin}_\ast^\text{op}\cong\text{Fin}^\text{op}\otimes\text{Fin}^\text{inj}$ -- (arbitrary,injective) morphisms;
\item $\text{Fin}^\text{iso}$ -- (bijective,bijective) morphisms;
\item $\text{Fin}_\ast^\text{inj}\cong\text{Fin}^\text{inj}\otimes\text{Fin}^{\text{inj},\text{op}}$ -- (injective,injective) morphisms (conditional on Conjecture \ref{ConjEnv}; see Example \ref{subBurn2});
\item $\text{Burn}^\text{eff}\cong\text{Fin}\otimes\text{Fin}^\text{op}$ -- (arbitrary,arbitrary) morphisms.
\end{itemize}
\end{remark}

As usual, the cofree functors play an important role.

\begin{remark}
The cofree functor $\text{Hom}(\text{Fin}_\ast^\amalg,-)$ is Segal's $\Gamma$-object construction \cite{SegalGamma}.

As in the proof of Proposition \ref{PropInj} (because $\text{Fin}^\text{inj}\cong\text{Env}(\mathbb{E}_0)$), the cofree functor $\text{Hom}(\text{Fin}^{\text{inj},\amalg},-)$ takes $\mathcal{C}^\odot$ to $\mathcal{C}_{1/}^\odot$.
\end{remark}

The free functor $\text{Fin}_\ast^\amalg\otimes -$ can be described by means of Theorem \ref{TensorFin}. That is, since $\text{Fin}_\ast^\amalg\otimes\mathcal{C}^\odot\cong\text{Fin}^\amalg\otimes\text{Fin}^{\text{inj},\text{op},\amalg}\otimes\mathcal{C}^\odot$, Theorem \ref{TensorFin} tells us that $\text{Fin}_\ast^\amalg\otimes\mathcal{C}^\odot$ is a full subcategory of $$\text{Hom}(\mathcal{C}^{\text{op},\odot}\otimes\text{Fin}^{\text{inj},\amalg},\text{Top}^\times)\cong\text{Hom}(\mathcal{C}^{\text{op},\odot},\text{Hom}(\text{Fin}^{\text{inj},\amalg},\text{Top}^\times))\cong\text{Hom}(\mathcal{C}^{\text{op},\odot},\text{Top}_\ast^\times),$$ where $\text{Top}_\ast$ is the $\infty$-category of pointed spaces.

The free functor $\text{Fin}^{\text{inj},\amalg}\otimes -$ is more mysterious, but see Conjecture \ref{ConjEnv}.

\subsection{Operads}
Our main examples of modules over $\text{Fin}^\text{inj}$ and $\text{Fin}_\ast$ arise from $\infty$-operads. For the rest of the section, we will look at these examples in more detail.

We will want to be careful with our notation. A (colored) $\infty$-operad $\mathcal{O}$ includes by definition the data of an $\infty$-category $[\mathcal{O}^\otimes]$ fibered over $\text{Fin}_\ast$. The fiber over the pointed set $\langle 1\rangle$ (a singleton plus a basepoint) is the $\infty$-category $\mathcal{O}$ of colors. More generally the objects of $[\mathcal{O}^\otimes]$ are labeled by tuples of objects of $\mathcal{O}$, and the fibration $[\mathcal{O}^\otimes]\rightarrow\text{Fin}_\ast$ sends an $n$-tuple to $\langle n\rangle$.

$[\mathcal{O}^\otimes]$ itself admits a symmetric monoidal structure $\odot$ given by agglutination of tuples. For example, $(A,B)\odot(C)\cong(A,B,C)$, and in particular, the fibration $[\mathcal{O}^\otimes]\rightarrow\text{Fin}_\ast$ lifts to a symmetric monoidal functor $[\mathcal{O}^\otimes]^\odot\rightarrow\text{Fin}_\ast^\amalg$.

This symmetric monoidal structure on $[\mathcal{O}^\otimes]$ is known to Lurie (\cite{HA} 2.2.4.7), but no details are given. We will give a conjectural universal property of $[\mathcal{O}^\otimes]^\odot$ (Conjecture \ref{ConjEnv}), but it is not important for us otherwise. A proof of this conjecture would likely also make clear how to construct $[\mathcal{O}^\otimes]^\odot$ explicitly.

\begin{example}
The empty tuple $()$ is a terminal object in $[\mathcal{O}^\otimes]$. That is, $[\mathcal{O}^\otimes]^\odot$ itself is a $\text{Fin}^{\text{inj},\text{op}}$-module, which is colored (in the sense of \ref{ColoredModule}).

In other words, we can think of $[\mathcal{O}^\otimes]^\odot$ as an `algebraic theory' relative to $\text{Fin}^{\text{inj},\text{op}}$. We will conjecture in \ref{ConjEnv} that it is precisely the $\text{Fin}^{\text{inj},\text{op}}$-algebraic theory which models $\mathcal{O}$-algebras.
\end{example}

\begin{example}
\label{OperadInj}
An $\infty$-operad $\mathcal{O}$ is called \emph{unital} if the empty tuple $()\in[\mathcal{O}^\otimes]$ is initial (\cite{HA} 2.3.1), or equivalently either of the following:
\begin{itemize}
\item $\text{Env}(\mathcal{O})^\odot$ is a $\text{Fin}^\text{inj}$-module;
\item $[\mathcal{O}^\otimes]^\odot$ is a $\text{Fin}^\text{inj}\otimes\text{Fin}^{\text{inj},\text{op}}$-module.
\end{itemize}
$\mathcal{O}$ is called \emph{reduced} if it is unital and the $\infty$-category $\mathcal{O}$ of colors is contractible. In this case, $\text{Env}(\mathcal{O})^\odot$ is a \emph{cyclic} $\text{Fin}^\text{inj}$-module.
\end{example}

The following proposition is a restatement into our algebraic language of \cite{HA} 2.4.3.9.

\begin{proposition}
\label{OpFin}
Suppose that $\mathcal{O}$ is an $\infty$-operad. Then there is an $\infty$-category $\mathcal{O}^\prime$ such that $$\text{Env}(\mathcal{O})^\odot\otimes\text{Fin}^\amalg\cong\text{Fin}[\mathcal{O}^\prime].$$ Here $\text{Fin}[\mathcal{O}^\prime]$ denotes the commutative Fin-algebra freely generated by the $\infty$-category $\mathcal{O}^\prime$.

If $\mathcal{O}$ is unital, then $\mathcal{O}^\prime\cong\mathcal{O}$, the underlying $\infty$-category of colors.
\end{proposition}

\begin{proof}
First, suppose that $\mathcal{O}$ is unital. It will suffice to show that the two sides are equivalent in $\text{CocartMon}_\infty$; that is, for every cocartesian monoidal $\infty$-category $\mathcal{C}^\odot$, $\text{Hom}(\text{Env}(\mathcal{O})^\odot\otimes\text{Fin}^\amalg,\mathcal{C}^\odot)\cong\text{Hom}(\text{Fin}[\mathcal{O}],\mathcal{C}^\odot)$. Indeed, $$\text{Hom}(\text{Env}(\mathcal{O})^\odot\otimes\text{Fin}^\amalg,\mathcal{C}^\odot)\cong\text{Hom}(\text{Env}(\mathcal{O})^\odot,\text{Hom}(\text{Fin}^\amalg,\mathcal{C}^\odot)),$$ which is equivalent by \cite{HA} 2.4.3.9 to $$\text{Hom}(\text{Env}(\mathcal{O})^\odot,\mathcal{C}^\odot)\cong\text{Alg}_\mathcal{O}(\mathcal{C}^\odot)\cong\text{Fun}(\mathcal{O},\mathcal{C})\cong\text{Hom}(\text{Fin}[\mathcal{O}],\mathcal{C}^\odot).$$ If $\mathcal{O}$ is not unital, note that $\mathcal{O}\otimes\mathbb{E}_0$ is a unital $\infty$-operad with symmetric monoidal envelope $\text{Env}(\mathcal{O})^\odot\otimes\text{Fin}^{\text{inj},\amalg}$, so $$\text{Env}(\mathcal{O})^\odot\otimes\text{Fin}^\amalg\cong\text{Env}(\mathcal{O}\otimes\mathbb{E}_0)^\odot\otimes\text{Fin}^\amalg\cong\text{Fin}[\mathcal{O}^\prime],$$ where $\mathcal{O}^\prime$ is the underlying $\infty$-category of $\mathcal{O}\otimes\mathbb{E}_0$.
\end{proof}

The author does not have an explicit description of $\text{Fin}[\mathcal{O}^\prime]$, but in general it should be much simpler than $\text{Env}(\mathcal{O})$ because it does not remember any of the operad structure. In the most important case, $\mathcal{O}^\prime$ is contractible and $\text{Fin}[\mathcal{O}^\prime]\cong\text{Fin}$. That is:

\begin{example}
The $\infty$-operad $\mathcal{O}$ is reduced if and only if $\text{Env}(\mathcal{O})^\odot$ is a cyclic $\text{Fin}^\text{inj}$-module such that $\text{Env}(\mathcal{O})^\odot\otimes\text{Fin}^\amalg\cong\text{Fin}$.
\end{example}

In Section 5.2, we will address to what extent we should expect a converse. That is, given a cyclic $\text{Fin}^\text{inj}$-module $\mathcal{M}$ with $\mathcal{M}\otimes\text{Fin}^\amalg\cong\text{Fin}$, does $\mathcal{M}$ arise from an $\infty$-operad? We conjecture that the answer to a slightly modified version of this question is yes.

We have seen that every operad $\mathcal{O}$ gives rise to a PROP $\text{Env}(\mathcal{O})^\odot$ and a Lawvere theory $\text{Env}(\mathcal{O})^\odot\otimes\text{Fin}^{\text{op},\amalg}$. However, there are certainly examples of Lawvere theories which do not arise from operads. We give the following example in terms of 1-categories and sets, but we expect a similar story for $\infty$-categories.

\begin{example}
The Lawvere theory $\mathcal{T}_R$ modeling $R$-modules \cite{Lawvere} does not arise from an operad (even colored nonunital) for any commutative ring $R$.

Indeed, suppose $\text{Mod}_R=\text{Alg}_{\mathcal{O}}(\text{Set}^\times)$, where $\mathcal{O}$ is an operad (possibly colored). We may assume $\mathcal{O}$ is unital, because $\mathcal{O}\otimes\mathbb{E}_0$ is a unital operad with the same $\text{Set}^\times$-algebras, thus the same Lawvere theory. Since $\mathcal{T}_R$ is semiadditive, $$\mathcal{T}_R\cong\mathcal{T}_R\otimes\text{Fin}^\amalg\cong(\text{Env}(\mathcal{O})^\odot\otimes\text{Fin}^{\text{op},\amalg})\otimes\text{Fin}^\amalg,$$ which is $\text{Burn}^\text{eff}[\mathcal{O}]$ by Proposition \ref{OpFin}. Thus $$\text{Mod}_R\cong\text{Hom}(\mathcal{T}_R,\text{Set}^\times)\cong\text{Fun}(\mathcal{O},\text{Hom}(\text{Burn}^{\text{eff},\amalg},\text{Set}^\times))\cong\text{Fun}(\mathcal{O},\text{CMon}).$$ However, $\text{Fun}(\mathcal{O},\text{CMon})$ cannot be additive (only semiadditive). To see this, note that there is a full subcategory inclusion $\text{CMon}\rightarrow\text{Fun}(\mathcal{O},\text{CMon})$ given by taking the constant functor. This inclusion is compatible with the semiadditive structure, but CMon is not additive.
\end{example}

Now suppose $\mathcal{O}$ is any $\infty$-operad (not necessarily unital). We call a morphism of $[\mathcal{O}^\otimes]$ \emph{inert} if it is of the form $X\odot f\colon X\odot Y\rightarrow X$ (up to equivalence), where $X$ and $Y$ are tuples and $f\colon Y\rightarrow()$ is a terminal map. On the other hand, call a morphism $f\colon S\rightarrow T$ of $\text{Fin}_\ast$ \emph{active} if $f^{-1}(\ast)=\{\ast\}$, and a morphism in $[\mathcal{O}^\otimes]$ \emph{active} if it lies over an active morphism in $\text{Fin}_\ast$.

\begin{remark}
Any morphism in $[\mathcal{O}^\otimes]$ can be factored uniquely as an inert morphism followed by an active morphism. That is, inert and active morphisms form a \emph{factorization system} for $[\mathcal{O}^\otimes]$ (\cite{HA} 2.1.2.4).

The subcategory $[\mathcal{O}^\otimes_\text{act}]$ (of all objects and only active morphisms) inherits a symmetric monoidal structure from $[\mathcal{O}^\otimes]$. This is the symmetric monoidal envelope $\text{Env}(\mathcal{O})^\odot\cong[\mathcal{O}^\otimes_\text{act}]^\odot$ (\cite{HA} 2.2.4).
\end{remark}

In a sense, $[\mathcal{O}^\otimes]$ is built from $\text{Env}(\mathcal{O})$ by `forcing the empty tuple $()$ to be terminal' -- that is, by forcing $\text{Env}(\mathcal{O})^\odot$ to be a $\text{Fin}^{\text{inj},\text{op}}$-module. Motivated by this observation, we make a conjecture.

\begin{conjecture}
\label{ConjEnv}
For an $\infty$-operad $\mathcal{O}$, there is an equivalence of symmetric monoidal $\infty$-categories $$[\mathcal{O}^\otimes]^\odot\cong\text{Env}(\mathcal{O})^\odot\otimes\text{Fin}^{\text{inj},\text{op},\amalg}.$$ Equivalently, $\text{Hom}([\mathcal{O}^\otimes]^\odot,\mathcal{C}^\odot)\cong\text{Alg}_\mathcal{O}(\mathcal{C}^\odot_{/1})$.
\end{conjecture}

\begin{example}
The $\mathbb{E}_\infty$-operad is given by $[\mathcal{O}^\otimes]=\text{Fin}_\ast$. The conjecture correctly asserts $\text{Fin}_\ast\cong\text{Fin}\otimes\text{Fin}^{\text{inj},\text{op}}$.
\end{example}

\begin{example}
\label{subBurn2}
The $\mathbb{E}_0$-operad is given by $[\mathcal{O}^\otimes]=\text{Fin}_\ast^\text{inj}$, the subcategory of $\text{Fin}_\ast$ given by all objects and only those morphisms $f\colon S\rightarrow T$ which are injective away from the basepoint (that is, $f(x)=f(y)\neq\ast$ implies $x=y$). Assuming the conjecture, $$\text{Fin}^\text{inj}\otimes\text{Fin}^{\text{inj},\text{op}}\cong\text{Fin}_\ast^\text{inj}.$$ This is consistent with Remark \ref{subBurn}: $\text{Fin}_\ast^\text{inj}$ can be identified with the subcategory of $\text{Burn}^\text{eff}$ spanned by all objects and only those morphisms such that both the ingressive and the egressive are injective.
\end{example}

\section{Connective spectra}
In the last section, we saw that various properties of symmetric monoidal $\infty$-categories (en route to full-fledged additive $\infty$-categories) correspond to the structure of modules over a solid semiring $\infty$-category.

In this section, we consider another example of a solid semiring $\infty$-category with dramatically different properties.

Since any Kan complex is a quasicategory, there is a full subcategory inclusion $\text{Top}\subseteq\text{Cat}_\infty$. This inclusion preserves products and has a right adjoint, so by Lemma \ref{EasyLemma}, there are full subcategory inclusions $$\text{Ab}_\infty\subseteq\text{CMon}_\infty\subseteq\text{SymMon}_\infty,$$ which are even symmetric monoidal. Recall from 2.1 that we are identifying $\text{Ab}_\infty$ with connective spectra. So connective spectra are examples of symmetric monoidal $\infty$-categories, and tensor products of the latter agree with smash products of the former.

If $E$ is a connective spectrum, we will write $\vec{E}$ to denote the symmetric monoidal $\infty$-category (emphasizing that we are thinking of paths in the $\infty$-groupoid $\Omega^\infty E$ as \emph{directed} morphisms; that is, in an $\infty$-category). Somewhat abusively, we say that such a symmetric monoidal $\infty$-category $\vec{E}$ `is' a connective spectrum.

Thus we may think of the sphere spectrum as a commutative semiring $\infty$-category $\vec{\mathbb{S}}$. In 4.1, we show that $\vec{\mathbb{S}}$-modules are exactly connective spectra; they recover the subcategory $\text{Ab}_\infty\subseteq\text{SymMon}_\infty$.

In other words, a symmetric monoidal $\infty$-category is an $\vec{\mathbb{S}}$-module if and only if it is \emph{both} grouplike (in the sense that every object has an inverse up to equivalence) and an $\infty$-groupoid (every morphism has an inverse up to equivalence).

Therefore the inclusion $\text{Ab}_\infty\rightarrow\text{SymMon}_\infty$ has both a left and a right adjoint, given by $\vec{\mathbb{S}}\otimes -$ and $\text{Hom}(\vec{\mathbb{S}},-)$, respectively. These are operations which produce spectra from symmetric monoidal $\infty$-categories; the former is related to algebraic K-theory and the latter to the Picard group. Thus, we may view the higher algebra of spectra as a special case of the commutative algebra of symmetric monoidal $\infty$-categories, with a variant of K-theory playing the role of tensoring up $\vec{\mathbb{S}}\otimes -$. This is the content of 4.2.

\subsection{The semiring category $\vec{\mathbb{S}}$}
As in the introduction to this section, connective spectra are examples of symmetric monoidal $\infty$-categories, and connective commutative ring spectra are examples of commutative semiring $\infty$-categories, via the symmetric monoidal full subcategory inclusion $$\text{Ab}_\infty\subseteq\text{CMon}_\infty\subseteq\text{SymMon}_\infty.$$ We will shortly see that this full subcategory $\text{Ab}_\infty$ can be recovered as modules over a solid semiring $\infty$-category, which is just the sphere spectrum $\vec{\mathbb{S}}$.

Unlike in Section 3, this is not a situation where the 1-categorical and $\infty$-categorical stories are verbatim the same. We will want to compare the 1-categorical story for motivation, even though it is in many ways \emph{less} transparent than the $\infty$-categorical one (see Remarks \ref{hS} and \ref{Z}).

\begin{example}
If $M$ is a commutative monoid (for example, an abelian group), then $M$ is a (discrete) symmetric monoidal category.

If $R$ is a commutative ring (or a commutative semiring), then $R$ is a (discrete) semiring category.

If $\mathcal{C}\in\text{SymMon}$ is in the full subcategory $\text{Ab}\subseteq\text{SymMon}$, we might say abusively that $\mathcal{C}$ \emph{is} an abelian group.
\end{example}

\begin{remark}
\label{SSolid}
If $R$ is a solid ring spectrum, $\vec{R}$ is also solid, since solidness is equivalent to idempotence under $\otimes$. Recall from Remark \ref{SolidRingSp} that the solid ring spectra have been classified (and are all connective): they are just the Moore spectra with $\pi_0 E$ isomorphic to a subring of $\mathbb{Q}$.

In particular, $\vec{\mathbb{S}}$ is a solid semiring $\infty$-category.
\end{remark}

\begin{warning}
\label{warning}
In contrast, solid rings are not necessarily solid as semiring categories. For example $\mathbb{Z}/2$ is a solid ring, but the map of semiring categories $\mathbb{Z}/2\otimes\mathbb{Z}/2\rightarrow\mathbb{Z}/2$ is \emph{not} an equivalence ($\pi_1(H\mathbb{Z}/2\wedge H\mathbb{Z}/2)\cong\mathbb{Z}/2$, not $0$), and therefore $\mathbb{Z}/2$ is not solid as a semiring category.

The difference is that the inclusion $\text{Top}\rightarrow\text{Cat}_\infty$ has a right adjoint, so Lemma \ref{EasyLemma} applies to show that $\text{Ab}_\infty\rightarrow\text{SymMon}_\infty$ is a symmetric monoidal functor. But $\text{Set}\rightarrow\text{Cat}$ does \emph{not} have a right adjoint.
\end{warning}

\begin{remark}
A symmetric monoidal $\infty$-category $\mathcal{C}$ is a connective spectrum if and only if both: all objects are invertible up to equivalence ($\mathcal{C}$ is grouplike); \emph{and} all morphisms are invertible up to equivalence ($\mathcal{C}$ is an $\infty$-groupoid).

But both of these conditions are simultaneously necessary: $\text{Fin}^\text{iso}$ is a symmetric monoidal $\infty$-category which is an $\infty$-groupoid but not grouplike. On the other hand, for an example of a symmetric monoidal category which is grouplike but not an $\infty$-groupoid, take the full subcategory of $\text{Mod}_R$ spanned by invertible modules (under $\otimes$).

If, however, $\mathcal{C}$ is a commutative semiring $\infty$-category which is grouplike (a \emph{ring $\infty$-category}), then $\mathcal{C}$ is an $\infty$-groupoid by Corollary \ref{CorGpoid} below.
\end{remark}

\begin{theorem}
\label{main3}
The sphere spectrum $\vec{\mathbb{S}}$ is a solid semiring $\infty$-category, and a symmetric monoidal $\infty$-category $\mathcal{C}$ admits an $\vec{\mathbb{S}}$-module structure if and only if $\mathcal{C}$ is a connective spectrum.
\end{theorem}

\begin{remark}
\label{hS}
We might ask whether a similar result holds for 1-categories. The answer is partially yes. However, we have already seen in Warning \ref{warning} that the functor $\text{Ab}\rightarrow\text{Cat}$ does not even have a right adjoint, so certainly Ab is not a category of modules in Cat.

Instead, we can only expect to recover `2-abelian groups' as a category of modules. These are symmetric monoidal categories which are both grouplike and groupoids.

The associated semiring category (which is itself an example of such a 2-abelian group) is $h\vec{\mathbb{S}}$, the homotopy category of $\vec{\mathbb{S}}$. As a groupoid, $h\vec{\mathbb{S}}$ has objects labeled by $\mathbb{Z}$, and all automorphism groups are isomorphic to $\pi_1\mathbb{S}\cong\mathbb{Z}/2$.

Since $\vec{\mathbb{S}}$ is a solid semiring $\infty$-category, $h\vec{\mathbb{S}}$ is a solid semiring 1-category, and a symmetric monoidal category admits an $h\vec{\mathbb{S}}$-module structure precisely if it is a groupoid with every object invertible up to isomorphism (a 2-abelian group). The proof is essentially the same as the proof of Theorem \ref{main3} (to follow).
\end{remark}

Now we are ready to prove Theorem \ref{main3}.

\begin{lemma}
\label{gpoid}
Every $\vec{\mathbb{S}}$-module is an $\infty$-groupoid.
\end{lemma}

\begin{proof}
Suppose $\mathcal{M}^{+}$ is an $\vec{\mathbb{S}}^{+,\times}$-module. As in the discussion before Example \ref{2.9}, $h\colon\text{SymMon}_\infty^\otimes\rightarrow\text{SymMon}^\otimes$ is symmetric monoidal, so $h\mathcal{M}^{+}$ is an $h\vec{\mathbb{S}}$-module (as 1-categories). It suffices to show that $h\mathcal{M}$ is a groupoid.

Let 0 and 1 denote the additive and multiplicative units of $h\vec{\mathbb{S}}$, and $-1$ an additive inverse of 1, with a chosen isomorphism $\alpha\colon(-1)+1\xrightarrow{\sim}0$. An integer $n$ denotes (as an object of $h\vec{\mathbb{S}}$) $1^{+n}$ if $n$ is positive or $(-1)^{+|n|}$ if $n$ is negative.

The module structure on $h\mathcal{M}$ induces symmetric monoidal functors $h\vec{\mathbb{S}}^{+}\otimes h\mathcal{M}^{+}\rightarrow h\mathcal{M}^{+}$ and therefore $$m_{-}\colon h\vec{\mathbb{S}}^{+}\rightarrow\text{Hom}(h\mathcal{M}^{+},h\mathcal{M}^{+})$$ such that $m_1$ is the identity functor and $m_0$ is the constant functor sending all of $h\mathcal{M}$ to the unit $0\in h\mathcal{M}$. Denote $m_{-1}X$ by $-X$, and note that $\alpha\colon(-1)+1\xrightarrow{\sim}0$ induces a natural isomorphism $\alpha_X\colon(-X)+X\xrightarrow{\sim}0$.

Suppose $f\colon X\rightarrow Y$ in $h\mathcal{M}$. The inverse to $f$ will be $-f\colon -Y\rightarrow -X$ after an appropriate shift. Specifically, the inverse to $f+0\colon X+0\rightarrow Y+0$ is the composition $$Y+0\xrightarrow{\alpha_X^{-1}+Y}X+(-X)+Y\xrightarrow{X+(-f)+Y}X+(-Y)+Y\xrightarrow{X+\alpha_Y}X+0.$$ To prove it is inverse amounts to a diagram chase; the diagram proving that $f^{-1}\circ f=\text{id}$ is as follows (the right square commutes by naturality of $\alpha_X$): $$\xymatrix{
X+0\ar[rrr]^-{\alpha_X^{-1}+X}\ar[d]_-{f} &&&X+(-X)+X\ar[r]^-{X+\alpha_X}\ar[d]_-{X+(-f)+f} &X+0\ar[d]^-{\text{id}} \\
Y+0\ar[r]_-{\alpha_X^{-1}+Y} &X+(-X)+Y\ar[rr]_-{X+(-f)+Y} &&X+(-Y)+Y\ar[r]_-{X+\alpha_Y} &X+0.
}$$
\end{proof}

\begin{proof}[Proof of theorem.]
We already saw that $\vec{\mathbb{S}}$ is solid (Example \ref{SSolid}).

Now because $\text{Ab}_\infty^\otimes\subseteq\text{SymMon}_\infty^\otimes$ as symmetric monoidal $\infty$-categories, any connective spectrum is an $\vec{\mathbb{S}}$-module (since it is an $\mathbb{S}$-module in $\text{Ab}_\infty^\otimes$). Conversely, suppose the symmetric monoidal $\infty$-category $\mathcal{C}^\odot$ is an $\vec{\mathbb{S}}$-module. By the lemma, $\mathcal{C}$ is a (symmetric monoidal) $\infty$-groupoid; that is, an $\mathbb{E}_\infty$-space. Since $\pi_0\mathcal{C}$ is a $\pi_0\mathbb{S}\cong\mathbb{Z}$ module, $\mathcal{C}$ is even grouplike, and therefore a connective spectrum.
\end{proof}

\begin{corollary}
\label{CorGpoid}
Suppose $\mathcal{R}$ is a commutative ring $\infty$-category -- that is, a commutative semiring $\infty$-category which is grouplike (every object has an additive inverse up to equivalence). Then $\mathcal{R}$ is an $\infty$-groupoid.
\end{corollary}

\begin{proof}
If $\mathcal{R}$ is a commutative ring $\infty$-category, $\mathcal{R}^\text{iso}$ is a commutative ring $\infty$-groupoid (connective commutative ring spectrum), and by Lemma \ref{EasyLemma} the inclusion $\mathcal{R}^\text{iso}\rightarrow\mathcal{R}$ is a semiring functor\footnote{Specifically, the semiring functor is the counit of the symmetric monoidal adjunction $i\colon\text{CMon}_\infty^{\wedge}\rightleftarrows\text{SymMon}_\infty^\otimes\colon(-)^\text{iso}$ evaluated at $\mathcal{R}$.}. The composite semiring functor $$\vec{\mathbb{S}}\rightarrow\mathcal{R}^\text{iso}\rightarrow\mathcal{R},$$ exhibits $\mathcal{R}$ as an $\vec{\mathbb{S}}$-algebra, thus an $\vec{\mathbb{S}}$-module. So $\mathcal{R}$ is an $\infty$-groupoid.
\end{proof}

In Remark \ref{hS}, we noted that the 1-categorical analogue of $\vec{\mathbb{S}}$ is the groupoid with objects labeled by $\mathbb{Z}$, and all automorphism groups isomorphic to $\mathbb{Z}/2$. In particular, it is not (as we might expect) $\mathbb{Z}$. Regardless, $\mathbb{Z}$ has some interesting properties as a semiring category:

\begin{example}
\label{Z}
We think of the ring $\mathbb{Z}$ as a commutative semiring $\infty$-category in one of two (equivalent) ways: either as the nerve of the discrete semiring 1-category, or as the connective Eilenberg-Maclane spectrum $\vec{H\mathbb{Z}}$. To avoid confusion, we will use the notation $\vec{\mathbb{Z}}$.

Since $\vec{\mathbb{Z}}$ is an $\vec{\mathbb{S}}$-module, any module over $\vec{\mathbb{Z}}$ is a connective spectrum. That is, $$\text{Mod}_{\vec{\mathbb{Z}}}(\text{SymMon}_\infty^\otimes)\cong\text{Mod}_{H\mathbb{Z}}(\text{Ab}_\infty^\wedge),$$ which is the $\infty$-category of chain complexes of abelian groups, concentrated in nonnegative degrees \cite{Shipley}.

But $\vec{\mathbb{Z}}$ is not solid, since $H\mathbb{Z}\wedge H\mathbb{Z}$ is the integral dual Steenrod algebra, and \emph{not} $H\mathbb{Z}$. On the other hand, the map $\pi_\ast(H\mathbb{Z}\wedge H\mathbb{Z})\rightarrow\pi_\ast H\mathbb{Z}$ is an isomorphism for $\ast=0,1$ and only stops being an isomorphism at $\ast=2$ (where $\pi_2(H\mathbb{Z}\wedge H\mathbb{Z})\cong\mathbb{Z}$).

The upshot is this: that $\mathbb{Z}$ \emph{is} solid as a semiring 1-category, but not as a semiring $n$-category for any $n>1$.
\end{example}

\begin{theorem}
\label{Z2}
The semiring (1-)category $\mathbb{Z}$ is solid, and a symmetric monoidal (1-)category $\mathcal{C}^\odot$ admits the structure of a $\mathbb{Z}$-module if and only if $\mathcal{C}^\odot$ is a 2-abelian group (both grouplike and a groupoid) and the symmetry isomorphisms $\sigma\colon X\odot Y\rightarrow Y\odot X$ in $\mathcal{C}$ are all identity morphisms when $X=Y$.
\end{theorem}

\begin{example}
If $A$ is any abelian group, the category $BA$ (with one object and morphisms labeled by $A$) is a symmetric monoidal category and a $\mathbb{Z}$-module.
\end{example}

\begin{lemma}
Let $\mathbb{N}^{+}$ be the discrete symmetric monoidal (1-)category corresponding to the commutative monoid of nonnegative integers under addition. As symmetric monoidal (1-)categories, $\mathbb{Z}^{+}\cong\mathbb{N}^{+}\otimes h\vec{\mathbb{S}}^{+}$.
\end{lemma}

\begin{proof}[Proof of lemma.]
Denote $\mathcal{T}=\mathbb{N}\otimes h\vec{\mathbb{S}}$. There are semiring functors from $\mathbb{N}$ and $h\vec{\mathbb{S}}$ to $\mathcal{T}$ given by (for example) $$h\vec{\mathbb{S}}\cong\text{Fin}^{\text{iso}}\otimes h\vec{\mathbb{S}}\rightarrow\mathbb{N}\otimes h\vec{\mathbb{S}}\cong\mathcal{T}.$$ First, we will analyze the image of $h\mathbb{S}^{+}$ in $\mathcal{T}$.

Since the `isomorphism classes' functor $\pi_0\colon\text{Cat}\rightarrow\text{Set}$ is product-preserving, likewise $\pi_0\colon\text{SymMon}^\otimes\rightarrow\text{CMon}^\otimes$ is symmetric monoidal. So $$\pi_0(\mathcal{T})\cong\mathbb{N}\otimes\mathbb{Z}\cong\mathbb{Z}.$$ Recall that every object of $h\vec{\mathbb{S}}$ has exactly one nontrivial automorphism. For $2\in h\vec{\mathbb{S}}$, this is the symmetry isomorphism of the symmetric monoidal structure $\sigma\colon 1\oplus 1\rightarrow 1\oplus 1$. For any other $n\in h\vec{\mathbb{S}}$, it is $\sigma\oplus(n-2)$. However, since the symmetric monoidal functor $\mathbb{N}^{+}\rightarrow\mathcal{T}$ preserves the symmetry automorphism, all symmetry automorphisms in $\mathcal{T}$ are identities. Therefore, every morphism of $h\vec{\mathbb{S}}$ is sent to the identity in $\mathcal{T}$.

Let $\mathcal{T}^\prime$ be the (semiring) subcategory of $\mathcal{T}$ spanned by the image of $h\vec{\mathbb{S}}$. We have just seen that $\mathcal{T}^\prime\cong\mathbb{Z}$. In particular, there are semiring functors $h\vec{\mathbb{S}}\rightarrow\mathcal{T}^\prime$ and $\mathbb{N}\rightarrow\mathcal{T}^\prime$, and therefore an induced semiring functor from the coproduct $\mathcal{T}\rightarrow\mathcal{T}^\prime$, such that $\mathcal{T}\rightarrow\mathcal{T}^\prime\subseteq\mathcal{T}$ is equivalent to the identity on $\mathcal{T}$. So $\mathcal{T}\cong\mathcal{T}^\prime\cong\mathbb{Z}$, as desired.
\end{proof}

\begin{proof}[Proof of theorem.]
We already saw in Remark \ref{Z} that $\mathbb{Z}$ is solid. Since it is itself a 2-abelian group (an $h\vec{\mathbb{S}}$-module), all $\mathbb{Z}$-modules are 2-abelian groups. Now say that a 2-abelian group is \emph{strictly commutative} if $\sigma\colon X\odot X\rightarrow X\odot X$ is the identity for all $X$. As in the proofs of Theorems \ref{main1} and \ref{main3}, we need only show both
\begin{enumerate}
\item $\text{Hom}(\mathbb{Z}^{+},\mathcal{C}^\odot)$ is strictly commutative for all symmetric monoidal $\mathcal{C}$;
\item if $\mathcal{C}^\odot$ is a strictly commutative 2-abelian group, $\text{Hom}(\mathbb{Z}^{+},\mathcal{C}^\odot)\rightarrow\mathcal{C}$ is an equivalence of categories.
\end{enumerate}

If $\mathcal{C}^\odot$ is any symmetric monoidal category, and $F\colon\mathbb{Z}^{+}\rightarrow\mathcal{C}^\otimes$ a symmetric monoidal functor, then by the definition of a symmetric monoidal functor, the following diagram commutes $$\xymatrix{
F(n)\odot F(n)\ar[r]^\sigma\ar[d]_{\epsilon} &F(n)\odot F(n)\ar[d]^{\epsilon} \\
F(n+n)\ar[r]_{=} &F(n+n).
}$$ Here the bottom morphism is the identity and $\epsilon$ is an isomorphism, so the symmetry isomorphism $\sigma$ in $\mathcal{C}^\odot$ is also the identity. Since this holds for every $n$, the symmetry map $\sigma\colon F\odot F\rightarrow F\odot F$ in the symmetric monoidal category $\text{Hom}(\mathbb{Z}^{+},\mathcal{C}^\odot)$ is also the identity. Therefore, $\text{Hom}(\mathbb{Z}^{+},\mathcal{C}^\odot)$ is strictly commutative.

For (2), suppose that $\mathcal{C}^\odot$ is a strictly commutative 2-abelian group. Note that being a strictly commutative 2-abelian group is invariant under equivalence of symmetric monoidal categories, so we may as well assume $\mathcal{C}^\odot$ is also permutative (or if not, replace it by a permutative category). Then $$\text{Hom}(\mathbb{Z}^{+},\mathcal{C}^\odot)\cong\text{Hom}(\mathbb{N}^{+},\text{Hom}(h\vec{\mathbb{S}}^{+},\mathcal{C}^\odot))\cong\text{Hom}(\mathbb{N}^{+},\mathcal{C}^\odot)$$ by the lemma and the fact that $\mathcal{C}^\odot$ is an $h\vec{\mathbb{S}}$-module (a 2-abelian group). So we need only show that the evaluation at 1 functor $\text{ev}_1\colon\text{Hom}(\mathbb{N}^{+},\mathcal{C}^\odot)\rightarrow\mathcal{C}$ is an equivalence of categories.

Consider the functor $F\colon\mathcal{C}\rightarrow\text{Hom}(\mathbb{N}^{+},\mathcal{C}^\odot)$ given by $$F(X)(n)=X^{\odot n}$$ $$F(f)(n)=f^{\odot n}\colon X^{\odot n}\rightarrow Y^{\odot n}.$$ Then $\text{ev}_1F$ is the identity, so $\text{ev}_1$ is full and essentially surjective. Moreover, $\text{ev}_1$ is faithful, and therefore an equivalence of categories, as desired.
\end{proof}

\subsection{Free and cofree spectra}
As with Fin and $\text{Fin}^\text{op}$, the cofree spectrum on a symmetric monoidal $\infty$-category, $\text{Hom}(\vec{\mathbb{S}}^{+},-)\colon\text{SymMon}_\infty\rightarrow\text{Ab}_\infty$, is not difficult to describe.

\begin{proposition}
\label{cofree2}
If $\mathcal{C}^\odot$ is a symmetric monoidal $\infty$-category, $\text{Hom}(\vec{\mathbb{S}}^{+},\mathcal{C}^\odot)$ is the symmetric monoidal subcategory of $\mathcal{C}^\odot$ spanned by those objects which have inverses (up to equivalence) and morphisms which have inverses (also up to equivalence). This is sometimes known as the \emph{Picard $\infty$-groupoid} $\text{Pic}(\mathcal{C}^\odot)$.
\end{proposition}

\begin{proof}
Let $\mathcal{D}$ be the full subcategory of $\mathcal{C}$ spanned by invertible objects (invertible up to equivalence). Then $\text{Hom}(\vec{\mathbb{S}}^{+},\mathcal{D}^\odot)\rightarrow\text{Hom}(\vec{\mathbb{S}}^{+},\mathcal{C}^\odot)$ is an equivalence, since every object of $\vec{\mathbb{S}}$ is invertible (and therefore is sent by any symmetric monoidal functor to an invertible object of $\mathcal{C}$).

Now we know there is a subcategory inclusion $$i\colon\text{Hom}(\vec{\mathbb{S}}^{+},\mathcal{D}^{\text{iso},\odot})\subseteq\text{Hom}(\vec{\mathbb{S}}^{+},\mathcal{D}^\odot)$$ because $\mathcal{D}^\text{iso}$ is a subcategory of $\mathcal{D}$. (This follows from the definition of symmetric monoidal functors \cite{HA} 2.1.3.7.) Any object of $\text{Hom}(\vec{\mathbb{S}}^{+},\mathcal{D}^\odot)$ factors through $\mathcal{D}^\text{iso}$ (is in the image of $i$, up to equivalence) since every morphism in $\vec{\mathbb{S}}$ is invertible. Moreover, because $\text{Hom}(\vec{\mathbb{S}}^{+},\mathcal{D}^\odot)$ is an $\vec{\mathbb{S}}$-module (Remark \ref{HomUp}), it is an $\infty$-groupoid by Lemma \ref{gpoid}. Thus any morphism $\phi\colon F\rightarrow G$ in $\text{Hom}(\vec{\mathbb{S}}^{+},\mathcal{D}^\odot)$ is in the image of $i$ (up to equivalence), since for any object $x\in\vec{\mathbb{S}}$, $\phi_x\colon F(X)\rightarrow G(X)$ must be an equivalence. Therefore $i$ is an equivalence.

But $\vec{\mathbb{S}}$ and $\mathcal{D}^{\text{iso},\odot}$ are both $\vec{\mathbb{S}}$-modules; write $\mathcal{D}^\text{iso}\cong\vec{E}$ where $E$ is a connective spectrum. Because $\vec{\mathbb{S}}$ is solid, $$\text{Hom}(\vec{\mathbb{S}}^{+},\mathcal{C}^\odot)\cong\text{Hom}(\vec{\mathbb{S}}^{+},\mathcal{D}^{\text{iso},\cdot})\cong\text{Hom}_{\vec{\mathbb{S}}}(\vec{\mathbb{S}}^{+},\mathcal{D}^{\text{iso},\cdot})\cong\text{Hom}_{\vec{\mathbb{S}}}(\vec{\mathbb{S}},\vec{E}).$$ Moreover, because $\text{Mod}_{\vec{\mathbb{S}}}^\otimes\cong\text{Ab}_\infty^\wedge$ as symmetric monoidal $\infty$-categories, we have $\text{Hom}_{\vec{\mathbb{S}}}(\vec{\mathbb{S}},\vec{E})=\vec{X}$, where $X=\text{Hom}_\text{Sp}(\mathbb{S},E)\cong E$. Thus $$\text{Hom}(\vec{\mathbb{S}}^{+},\mathcal{C}^\odot)\cong\vec{E}\cong\mathcal{D}^\text{iso},$$ as desired.
\end{proof}

As in Section 3.2, the free construction $\vec{\mathbb{S}}\otimes -$ is much harder to compute, but we can say something nonetheless. In this case, we think of $\vec{\mathbb{S}}\otimes\mathcal{C}^\odot$ as a variant of algebraic K-theory of $\mathcal{C}^\odot$. Where the usual algebraic K-theory of a symmetric monoidal $\infty$-category roughly involves taking the core (throwing out noninvertible morphisms) and then the $\infty$-group completion, $\vec{\mathbb{S}}\otimes -$ involves taking the classifying space (formally inverting all morphisms) and then the $\infty$-group completion. Of course, if $\mathcal{C}^\odot$ is already an $\infty$-groupoid, the two constructions agree.

\begin{proposition}
If $\mathcal{C}^\odot$ is symmetric monoidal, then $\vec{\mathbb{S}}\otimes\mathcal{C}^\odot\cong K(|\mathcal{C}^\odot|)$. Here $|\mathcal{C}|$ is the geometric realization of the $\infty$-category $\mathcal{C}$ (construed as a quasicategory and thus a simplicial set), and $K(-)$ is the $\infty$-group completion of an $\mathbb{E}_\infty$-space.
\end{proposition}

\begin{proof}
By Theorem \ref{main3}, $\vec{\mathbb{S}}\otimes -$ is left adjoint to the subcategory inclusion $\text{Ab}_\infty\rightarrow\text{SymMon}_\infty$, which factors $$\text{Ab}_\infty\rightarrow\text{CMon}_\infty\rightarrow\text{SymMon}_\infty.$$ Therefore, $\vec{\mathbb{S}}\otimes -$ factors as a composition of left adjoints $$\xymatrix{
\text{SymMon}_\infty\ar[r]^{|-|} &\text{CMon}_\infty\ar[r]^{K(-)} &\text{Ab}_\infty.
}$$
\end{proof}

For example, if $\mathcal{C}^\odot$ has an initial or terminal object, its classifying space is contractible. The following corollary is an algebraic restatement of this classical fact.

\begin{corollary}
\label{vanish}
$\vec{\mathbb{S}}\otimes\text{Fin}^{\text{inj}}\cong 0\cong\vec{\mathbb{S}}\otimes\text{Fin}^{\text{inj},\text{op}}$. As a result, $\vec{\mathbb{S}}\otimes\mathcal{C}^\odot=0$ for any $\mathcal{C}$ whose unit is initial or terminal, and in particular for $\mathcal{C}^\odot$ any of the solid semiring $\infty$-categories discussed before: Fin, $\text{Fin}_\ast$, $\text{Fin}^\text{inj}$, $\text{Burn}^\text{eff}$, or their opposites.
\end{corollary}

We also have a slightly stronger result.

\begin{corollary}
\label{VanishCor}
If $\mathcal{C}^\odot$ is a $\text{Fin}^\text{inj}$-module, then the only object which is invertible (up to equivalence) is the unit.
\end{corollary}

\begin{proof}
By Proposition \ref{cofree2}, $\text{Hom}(\vec{\mathbb{S}}^{+}\otimes\text{Fin}^{\text{inj},\amalg},\mathcal{C}^\odot)\cong\text{Hom}(\vec{\mathbb{S}}^{+},\mathcal{C}^\odot)$ is the maximal subgroupoid spanned by invertible objects. But it is also a module over $\vec{\mathbb{S}}\otimes\text{Fin}^\text{inj}\cong 0$, and therefore contractible. So every invertible object is equivalent to the unit.
\end{proof}

\begin{example}
Let $\mathcal{H}^\odot$ be the following pushout of symmetric monoidal (not semiring!) $\infty$-categories (all symmetric monoidal structures are disjoint union, as usual) $$\xymatrix{
\text{Fin}^\text{iso}\ar[r]\ar[d] &\text{Fin}\ar[d] \\
\text{Fin}^\text{op}\ar[r] &\mathcal{H}.
}$$ In particular, an object of $\text{Hom}(\mathcal{H}^\odot,\mathcal{C}^\odot)$ consists of an object of $\mathcal{C}$ along with both a commutative algebra and a cocommutative coalgebra structure (but not necessarily compatible with each other). Tensoring with symmetric monoidal $\infty$-categories $\mathcal{C}^\odot\otimes -$ preserves colimits (such as pushouts) because the functor has a right adjoint $\text{Hom}(\mathcal{C}^\odot,-)$. Thus $\mathcal{H}^\odot\otimes\vec{\mathbb{S}}$ is the pushout $$\xymatrix{
\vec{\mathbb{S}}\ar[r]\ar[d] &0\ar[d] \\
0\ar[r] &\mathcal{H}^\odot\otimes\vec{\mathbb{S}}.
}$$ Since this is a square of $\vec{\mathbb{S}}$-modules (connective spectra), we can take the pushout in spectra, so that $\mathcal{H}^\odot\otimes\vec{\mathbb{S}}\cong\vec{\Sigma\mathbb{S}}$, the suspension of $\mathbb{S}$.

On the other hand, $\mathcal{H}^\odot\otimes\text{Fin}^\amalg\cong\mathcal{H}^\odot\otimes\text{Fin}^{\text{op},\amalg}\cong\text{Burn}^\text{eff}$, by analyzing the same pushout square.
\end{example}

\section{Questions and conjectures}
\subsection{Computation of tensor products}
A major obstacle in this subject is the difficulty in making any computations, such as the computation of tensor products of symmetric monoidal $\infty$-categories. In Sections 3.2 and 4.2, we saw strategies for computing $\mathcal{C}\otimes\text{Fin}$, $\mathcal{C}\otimes\text{Fin}^\text{op}$, and $\mathcal{C}\otimes\vec{\mathbb{S}}$, which are useful at least for some symmetric monoidal $\mathcal{C}$.

We also know in many examples that tensor products are related to the span (or effective Burnside) construction. Notably:
\begin{itemize}
\item $\text{Fin}\otimes\text{Fin}^\text{op}\cong\text{Burn}^\text{eff}$ (Corollary \ref{main2});
\item all the other semiring $\infty$-categories of Section 3 can be described as span subcategories of $\text{Burn}^\text{eff}$ (Remark \ref{subBurn} and, conjecturally, Example \ref{subBurn2})
\item there is some reason to believe $\text{Env}(\mathcal{O})\otimes\text{Fin}^{\text{op},\amalg}$ can be described via a span construction:
\end{itemize}

\begin{remark}
The last point requires justification. Let $\mathcal{O}$ be a reduced $\infty$-operad, and recall from Example \ref{OperadCompute} that the associated Lawvere theory has mapping spaces $$\text{Map}(X^{\otimes m},X^{\otimes n})\cong (\coprod_{T\in\text{Fin}}\mathcal{O}(T)\times_{\Sigma_T}m^T)^n.$$ This mapping space can equivalently be described as the space of spans $$\xymatrix{
&T\ar[ld]_f\ar[rd]^g &\\
[m] &&[n],
}$$ where $T$ is an arbitrary finite pointed set, $f$ is an active morphism of finite pointed sets, and $g$ is an active morphism in $[\mathcal{O}^\otimes]$. Thus we expect that the associated Lawvere theory can be described in terms of a span construction.
\end{remark}

A very natural question then is:

\begin{question}
Is there a general class of tensor products of symmetric monoidal $\infty$-categories which can be computed via span constructions?
\end{question}

On the other hand, we might approach this question from the opposite direction. Instead of trying to compute tensor products, we might regard tensor product formulas like $\text{Fin}\otimes\text{Fin}^\text{op}\cong\text{Burn}^\text{eff}$ as providing algebraic universal properties for span constructions. There are many places where span constructions naturally arise, and a universal property would be helpful. A particularly natural example comes from equivariant homotopy theory, as in Example \ref{EqEx}.

\begin{question}
Let $G$ be a finite group, $\text{Fin}_G$ the category of finite $G$-sets, and $\text{Burn}_G^\text{eff}$ the category of spans of finite $G$-sets. Can $\text{Burn}_G^\text{eff}$ be decomposed as a tensor product? Is it true that $$\text{Burn}_G^\text{eff}\cong\text{Fin}_G\otimes_{\text{Fin}_G^\text{iso}}\text{Fin}_G^\text{op}?$$
\end{question}

This last question is equivalent to asking whether $$\text{Hom}(\text{Burn}_G^{\text{eff},\amalg},\text{Top}^\times)\cong\text{Hom}_{\text{Fin}_G^\text{iso}}(\text{Fin}_G^\amalg,\text{Hom}(\text{Fin}_G^{\text{op},\amalg},\text{Top}^\times));$$ or whether equivariant $\mathbb{E}_\infty$-spaces coincide with `$G$-commutative monoids' in equivariant spaces. Here we mean $G$-commutative monoids in the sense of Mazur \cite{Mazur} and Kaledin \cite{Kaledin}. The more refined notions of Hill and Hopkins \cite{HillHopkins} or Barwick et al. \cite{Barwick} require that we work with \emph{equivariant} symmetric monoidal $\infty$-categories. In that fully equivariant setting, we very much expect to have a statement of the form $\text{Burn}_G^\text{eff}\cong\text{Fin}_G\otimes_{\text{Fin}_G^\text{iso}}\text{Fin}_G^\text{op}$, but the successes of \cite{Mazur} and \cite{Kaledin} suggest that there is some hope of proving such a statement even for ordinary symmetric monoidal $\infty$-categories.

\subsection{Reconstruction of operads}
We have seen that to every unital $\infty$-operad $\mathcal{O}$ (given by a fibration $[\mathcal{O}^\otimes]\rightarrow\text{Fin}_\ast$) is associated a PROP, or symmetric monoidal envelope, satisfying $\text{Env}(\mathcal{O})^\odot\otimes\text{Fin}^\amalg\cong\text{Fin}[\mathcal{O}]$ (Proposition \ref{OpFin}).

\begin{remark}
The associated PROP is colored (see Remark \ref{ColoredModule}); that is, it includes the data of a full subcategory $\mathcal{O}$ inducing an essentially surjective symmetric monoidal functor $\text{Fin}^\text{iso}[\mathcal{O}]\rightarrow\text{Env}(\mathcal{O})^\odot$. Proposition \ref{OpFin} asserts that this map becomes an equivalence after tensoring with $\text{Fin}$; we say that the colored PROP is \emph{trivial over Fin}.
\end{remark}

We might now ask whether we can reconstruct an $\infty$-operad from the associated colored PROP. Assuming Conjecture \ref{ConjEnv}, the answer is yes:

There is a symmetric monoidal functor $\text{Fin}[\mathcal{O}]\rightarrow\text{Fin}$ induced from the trivial functor $\mathcal{O}\rightarrow\text{Fin}$ sending all of $\mathcal{O}$ to the singleton set. Then we have a composite $$\text{Env}(\mathcal{O})\rightarrow\text{Env}(\mathcal{O})^\odot\otimes\text{Fin}^\amalg\cong\text{Fin}[\mathcal{O}]\rightarrow\text{Fin}.$$ Tensoring this functor with $\text{Fin}^{\text{inj},\text{op}}$ produces (assuming Conjecture \ref{ConjEnv}) a functor $[\mathcal{O}^\otimes]\rightarrow\text{Fin}_\ast$ which recovers the original $\infty$-operad.

This motivates a second question: given any colored PROP $\mathcal{T}$ which is trivial over Fin, will this procedure produce an $\infty$-operad? The answer is no, for the following reason: if $\mathcal{T}\otimes\text{Fin}^{\text{inj},\text{op}}$ were an $\infty$-operad, it should have the property that picking an object over $[2]$ amounts to picking a pair of objects, each over $[1]$ (condition (2) of Definition \cite{HA} 2.1.1.10). But there is no reason to suspect such a thing, \emph{unless} $\mathcal{T}$ is suitably built out of a cartesian monoidal $\infty$-category (a Lawvere theory):

\begin{conjecture}
The functor $\text{Env}(-)^\odot\otimes\text{Fin}^\text{op}$, from unital $\infty$-operads to colored $\text{Fin}_\ast^\text{op}$-modules (`colored pointed Lawvere theories', if you like), is fully faithful. A colored $\text{Fin}_\ast^\text{op}$-module is in the image of this functor if and only if it is trivial over $\text{Burn}^\text{eff}$.
\end{conjecture}

\section*{Acknowledgments} 
\thispagestyle{empty}
This paper has benefited from the comments and advice of a number of people, including Mike Hill, Clark Barwick, Saul Glasman, and Ben Knudsen, as well as many more. The author is also thankful to the referee for reading the paper quickly and carefully, with many helpful comments and suggestions.

\end{document}